\documentclass[11pt,a4paper]{article}
\usepackage[round,authoryear,comma]{natbib}
\usepackage{optimstyle}
\usepackage{tabularx,pifont,calc,fontawesome5}

\newcommand{\R}{\mathbb R}
\newcommand{\F}{\mathcal F}
\newcommand{\ip}[2]{\langle #1,#2\rangle}
\newcommand{\norm}[1]{\|#1\|}
\DeclareMathOperator{\argmin}{arg\,min}

\newcolumntype{P}[1]{>{\raggedright\arraybackslash}p{#1}}
\newcolumntype{Y}{>{\raggedright\arraybackslash}X}

\newcommand{\yes}{\textcolor{green!55!black}{\ding{51}}}
\newcommand{\no}{\textcolor{red!80!black}{\ding{55}}}

\setlist[itemize]{label=$\vartriangleright$}
\newcommand{\by}[2]{\mathrel{\underset{#2}{#1}}}
\newcommand{\bypad}[3]{\mathrel{\underset{\mathmakebox[\widthof{$\scriptstyle#3$}]{#2}}{#1}}}
\BeforeBeginEnvironment{theorem}{\Needspace{6\baselineskip}}
\BeforeBeginEnvironment{lemma}{\Needspace{6\baselineskip}}
\BeforeBeginEnvironment{proposition}{\Needspace{6\baselineskip}}
\BeforeBeginEnvironment{definition}{\Needspace{6\baselineskip}}

\renewcommand{\optimtitlefont}{\fontsize{16}{21}\selectfont\bfseries\boldmath}
\title{Descending Everest, Part I\texorpdfstring{\\[3pt]}{: } Oracle Complexity of Convex Optimization\texorpdfstring{\\}{ }under \texorpdfstring{$(L_0,L_1)$}{(L0,L1)}-Smoothness\texorpdfstring{\thanks{The Russian version of this paper is available at \url{https://github.com/chitoraphelion/descending-everest-with-first-order-oracles}.}}{}}
\author{%
  Mikhail Rudomanov\thanks{MIRIAI, Moscow Independent Research Institute of Artificial Intelligence, Moscow, Russia. \email{rudomanov.ma@phystech.edu}.}\,\thanks{Equal contribution.}\quad
  Ivan Fomin\thanks{Yandex, Moscow, Russia. \email{ash-ivanfomin@yandex-team.ru}.}\,\footnotemark[3]\quad
  Alexander V. Gasnikov\thanks{Innopolis University, Innopolis, Russia; \email{gasnikov@yandex.ru}.}\quad
  Fedor S. Stonyakin\thanks{V.\,I.~Vernadsky Crimean Federal University, Simferopol, Russia. \email{fedyor@mail.ru}.}}
\date{}

\begin{document}

\begingroup
\makeatletter
\renewcommand{\@fnsymbol}[1]{\ensuremath{\ifcase#1\or *\or\spadesuit\or\text{\tiny\faMountain}\or\clubsuit\or\dagger\or\ddagger\fi}}
\let\origmakefntext\@makefntext
\long\def\@makefntext#1{\rightskip\z@\@plus3em\relax\origmakefntext{#1}}
\makeatother
\maketitle
\endgroup

\begin{abstract}
For convex $(L_0,L_1)$-smooth functions, the curvature may grow together with the gradient norm. Therefore, even when $\norm{x^0-x^\ast}\le R$, the gradient and the optimality gap at the starting point may be exponentially large in $L_1R$, and the bound $\sqrt{L_0R^2/\varepsilon}$, known for the accelerated method on functions with Lipschitz gradient, no longer describes the complexity of the problem. We determine the oracle complexity of this class for arbitrary deterministic first-order methods: for $L_0R^2/\varepsilon\ge1$, the optimal number of oracle calls $N$ satisfies $1+N=\widetilde\Theta\bigl(\sqrt{L_0R^2/\varepsilon}+L_1R\bigr)$. The upper bound reads $O\bigl(\sqrt{L_0R^2/\varepsilon}+(1+L_1R)\log(2+L_1R)\bigr)$: the method first passes through the large-gradient regime in $O(1+L_1R)$ iterations with an inexact line search, and then the standard accelerated method with a constant stepsize takes over. The lower bound is obtained by means of chains and the orthogonal invariance principle on functions conjugate to strongly convex functions on a polytope. If the initial optimality gap is at most $4\varepsilon$, the optimal bound becomes $1+N=\widetilde\Theta\bigl(\sqrt{L_0R^2/\varepsilon}+(L_1R)^{2/3}\bigr)$.
\end{abstract}

\begin{tcolorbox}[
  optimbox,
  colback=optimBlue!3!white,
  colframe=optimBlue!40!white,
  borderline west={1.4pt}{0pt}{optimBlue},
  top=7pt,bottom=7pt
]
\small
\noindent\textbf{Use of language models.}\enspace
The models \textit{Codex Astra 6} and \textit{Claude Opus 5.5} were used
in obtaining the theoretical results of this work.
They were also used in preparing and typesetting the manuscript and in
translating the original text into English.
However, every line of this paper, including every statement, every
proof and every computation, has been personally edited and verified
by the authors. The authors bear sole responsibility for the correctness
of all mathematical statements, proofs and conclusions.
\end{tcolorbox}

\section{Introduction}\label{sec:intro}

We consider the unconstrained convex optimization problem
\[
 f(x)\to\min_{x\in\R^d},
\]
where $f:\R^d\to\R$ is a convex differentiable $(L_0,L_1)$-smooth function that attains its minimum $f^\ast$ at some point\footnote{The minimizer is, in general, not unique. The set $X^\ast=\argmin f$ is convex and closed, so the condition below is equivalent to $\operatorname{dist}(x^0,X^\ast)\le R$, and $x^\ast$ may be taken to be the minimizer closest to $x^0$.} $x^\ast$. About the starting point $x^0$, only an upper bound $R$ on its distance to a minimizer is known: $\norm{x^0-x^\ast}\le R$.

For convex problems with $L$-Lipschitz gradient, the distance $R$ to the minimizer bounds the initial optimality gap by $LR^2/2$ and determines the complexity of gradient methods, namely the classical bounds $LR^2/\varepsilon$ and $\sqrt{LR^2/\varepsilon}$. Under generalized smoothness this is no longer the case. For twice continuously differentiable functions, the class of $(L_0,L_1)$-smooth convex functions with parameters $L_0,L_1>0$ is defined by the condition
\[
 0\preceq\nabla^2f(x)\preceq\bigl(L_0+L_1\norm{\nabla f(x)}\bigr)I ,
\]
that is, the curvature is allowed to grow together with the gradient norm. This model was proposed to explain the success of gradient clipping methods \citep{zhang2019gradient,zhang2020improved} and has since become one of the standard assumptions in the analysis of first-order methods.

The $(L_0,L_1)$-smoothness condition was proposed on the basis of empirical observations in neural network training. Experiments with an LSTM language model revealed a positive correlation between the local smoothness constant and the gradient norm, which motivated the model with the dependence $L_0+L_1\norm{\nabla f(x)}$ \citep{zhang2019gradient}. A similar dependence was observed in transformer training \citep{crawshaw2022robustness}.

The $(L_0,L_1)$-smoothness condition does not require convexity. In particular, it holds for univariate polynomials of degree at least three with suitable constants $L_0,L_1$ \citep{zhang2019gradient}. Moreover, the corresponding Hessian bound holds for some distributionally robust optimization problems after passing to the dual variable and adding a penalty, under certain assumptions on the loss function and the divergence \citep{jin2021non,chen2023generalized}.

Unlike for functions with Lipschitz gradient, the distance to the minimizer no longer bounds the optimality gap by a quadratic quantity: at a point at distance $R$ from the minimizer, the gap may be as large as $(L_0/L_1^2)\bigl(e^{L_1R}-L_1R-1\bigr)$, that is, exponentially large in $L_1R$, and this bound is attained (Lemma~\ref{lem:radial}). Therefore, before the method reaches the region of small gradients, where the familiar theory for problems with Lipschitz gradient applies, it has to pass through the region of large gradients. How many oracle calls does this cost the best first-order method if only an upper bound $R$ on the distance from the starting point to a minimizer is known? In this paper we answer this question: up to logarithmic factors, the cost is $L_1R$.

\paragraph{Notation.} We denote by $\varepsilon>0$ the required accuracy in function value: the method must output a point $\hat x$ with $f(\hat x)-f^\ast\le\varepsilon$. The symbol $\norm{\cdot}$ always stands for the Euclidean norm of a vector and the spectral norm of a matrix, $[x,y]$ is the segment with endpoints $x,y\in\R^d$, $B(x,r)$ is the closed ball of radius $r$ centered at $x$, and C--S stands for the Cauchy--Schwarz inequality. The sum of the first $k$ terms of a sequence is denoted by the corresponding capital letter: for example, $A_k=a_1+\dots+a_k$ and $A_0=0$.

\paragraph{Main result.} We consider arbitrary deterministic first-order methods: a method may query the gradient at any points and output any point, and every oracle call is counted. We prove that the smallest number of calls $N$ for which there exists a deterministic first-order method solving the problem for all functions of the class in all dimensions satisfies\footnote{Here and below, the lower bounds are proved for $L_0R^2\ge\varepsilon$.}
\begin{align*}
 1+N&\ge\frac1{29}\Bigl[1+\sqrt{\frac{L_0R^2}{\varepsilon}}+\frac{L_1R}{\sqrt{\log(2+L_1R)}}\Bigr],\\
 1+N&\le5\sqrt{\frac{L_0R^2}{\varepsilon}}+C\,(1+L_1R)\log(2+L_1R)
\end{align*}
with an absolute constant $C$ (Theorems~\ref{thm:upper} and~\ref{thm:lower-radius}). Thus, the optimal complexity is $\sqrt{L_0R^2/\varepsilon}+L_1R$ up to logarithmic factors, and the logarithmic factors affect only the $L_1R$ term. The method achieving the upper bound uses only gradients. The lower bound remains valid for methods that receive the function value together with the gradient.

\paragraph{Upper and lower bounds.} The upper bound is obtained by a two-stage method. At the first stage, the method SHERPA passes through the large-gradient regime in $O(1+L_1R)$ iterations. This is an accelerated method with the inexact line search ICEPICK, based on Nesterov's framework of estimating sequences \citep{nesterov2018lectures}. At the second stage, we run the standard accelerated method with a constant stepsize $1/(2L_0)$, using the analysis of \citet{tyurin2026accelerated}. We call the variant of this method that we use CRUISE. The lower bound follows the classical scheme: a chain and the orthogonal invariance principle \citep{nemirovski1983problem,nesterov2018lectures,carmon2020lower}.

\paragraph{Bounded initial optimality gap.} If it is additionally known that $f(x^0)-f^\ast\le4\varepsilon$, the complexity drops to
\[
 1+N=\widetilde\Theta\Bigl(\sqrt{\frac{L_0R^2}{\varepsilon}}+(L_1R)^{2/3}\Bigr)
\]
(Appendix~\ref{app:gap}). There is no contradiction with the main result: the worst-case functions that yield the term $L_1R$ have a huge initial optimality gap and do not belong to the class with bounded gap. This regime may be useful for restart schemes: each subsequent run of the method starts from a point at which the optimality gap is already bounded by a constant multiple of the accuracy required from this run, so the cost of one run is governed by the exponent $2/3$ rather than $1$.

\paragraph{Related work.} The generalization in which $L_0+L_1\norm{\nabla f(x)}$ is replaced by an arbitrary function $\ell(\norm{\nabla f(x)})$ was studied by \citet{li2023generalized,tyurin2025unified}. We restrict ourselves to the affine case. Convergence guarantees for gradient methods in the convex case were obtained by \citet{gorbunov2024methods,vankov2025optimizing,lobanov2025linear}, and accelerated methods were considered by \citet{gorbunov2024methods,vankov2025optimizing,tyurin2026accelerated,borodich2026graal}. Lower bounds for functions with Lipschitz gradient and arbitrary first-order methods are discussed, for example, by \citet{guzman2015lower}. A summary table of known upper bounds for the convex case is given in Appendix~\ref{app:literature}.

\paragraph{Comparison with known results.} Known bounds for the convex deterministic case that do not use a bound on the initial optimality gap contain either a term $(L_1R)^2$ or higher, or a factor equal to the logarithm of the ratio of the initial gap $f(x^0)-f^\ast$ to $\varepsilon$ or to $L_0/L_1^2$ (Table~\ref{tab:lit-convex}). On the radius-only class, this factor may be of order $L_1R$. For example, when $L_1R\ge1$, to pass through the large-gradient regime, that is, to reach an optimality gap of at most $L_0/(128L_1^2)$, the known bounds require $O((L_1R)^{5/3})$ iterations with exact line search for the AGMsDR method \citep{vankov2025optimizing} and $O((L_1R)^2)$ iterations for the preliminary stage of the method of \citet{tyurin2026accelerated}. Our first stage passes through this regime in $O\bigl((1+L_1R)\log(2+L_1R)\bigr)$ oracle calls with unknown initial optimality gap. After that, the accelerated method of \citet[Theorem~3.2]{tyurin2026accelerated} gives the bound $O\bigl(1+\sqrt{L_0R^2/\varepsilon}\bigr)$ without logarithms, and we use it at the second stage. Lower bounds in the generalized setting were previously known only for specific methods, for example, for gradient descent with a constant stepsize \citep{zhang2019gradient,crawshaw2022robustness}, variants of AdaGrad \citep{crawshaw2025lower} and parameter-agnostic methods \citep{hubler2024parameter}, or were inherited from the subclass $L_1=0$ \citep{nesterov2018lectures}. To the best of our knowledge, Theorems~\ref{thm:lower-radius} and~\ref{thm:lower-gap} are the first lower bounds for arbitrary deterministic methods that involve $L_1$.

\paragraph{Organization.} Section~\ref{sec:setting} defines the function class and the oracle complexity. Section~\ref{sec:upper} constructs the method that gives the upper bound. Section~\ref{sec:lower} proves the lower bound, which matches the upper one asymptotically up to logarithmic factors. Section~\ref{sec:conclusion} lists open questions. Appendix~\ref{app:tools} contains general inequalities for functions of the class described above and the proofs of auxiliary lemmas. The appendices also contain the class with bounded initial optimality gap (Appendix~\ref{app:gap}), fixed dimension (Appendix~\ref{app:fixed-dim}) and a summary table of known bounds (Appendix~\ref{app:literature}).

\begin{table}[ht]
\centering\small
\caption{Results obtained in this paper. The column ``Condition'' gives the additional restriction on the initial optimality gap in the definition of the class.}
\label{tab:results}
\begin{tabularx}{\linewidth}{@{}P{.2\linewidth}lccY@{}}
\toprule
Result & Condition & Lower & Upper & Complexity $1+N$\\
\midrule
Theorems~\ref{thm:upper}, \ref{thm:lower-radius} & none & \yes & \yes & $\widetilde\Theta\Bigl(1+\sqrt{\dfrac{L_0R^2}\varepsilon}+L_1R\Bigr)$\\[3mm]
Theorems~\ref{thm:upper-gap}(a), \ref{thm:lower-gap} & $f(x^0)-f^\ast\le4\varepsilon$ & \yes & \yes & $\widetilde\Theta\Bigl(1+\sqrt{\dfrac{L_0R^2}\varepsilon}+(L_1R)^{2/3}\Bigr)$\\[3mm]
Theorem~\ref{thm:upper-gap}(b) & $f(x^0)-f^\ast\le\Delta$ & \no & \yes & $O\Bigl(\log\Bigl(2+\dfrac{L_0R^2}\varepsilon+L_1R\Bigr)\,\times$\newline$\times\Bigl[1+\sqrt{\dfrac{L_0R^2}\varepsilon}+(L_1R)^{2/3}\log^{2/3}\dfrac{2\Delta}\varepsilon\Bigr]\Bigr)$\\
\bottomrule
\end{tabularx}
\par\smallskip\raggedright\footnotesize
The third row assumes $\Delta>\varepsilon$. For $\Delta\le\varepsilon$ it suffices to return $x^0$ without oracle calls.
\end{table}

\section{Problem setting}\label{sec:setting}

In this section we introduce the function class and the oracle complexity.

\begin{definition}[Function class]\label{def:class}
Let $\ell:[0,\infty)\to(0,\infty)$ be a continuous nondecreasing function. We denote by $\F_\ell^d$ the class of convex differentiable functions $f:\R^d\to\R$ such that for all $x,y\in\R^d$
\begin{equation}\label{eq:class}
 \norm{\nabla f(x)-\nabla f(y)}
 \le\ell\Bigl(\sup_{u\in[x,y]}\norm{\nabla f(u)}\Bigr)\norm{x-y} .
\end{equation}
For $\ell(s)=L_0+L_1s$ we write $\F_{L_0,L_1}^d$ and call the functions of this class \emph{$(L_0,L_1)$-smooth}.
\end{definition}

Inequality~\eqref{eq:class} is the symmetric $(L_0,L_1)$-smoothness condition of \citet[Assumption~1.3]{gorbunov2024methods}, written for an arbitrary function $\ell$. It agrees with the definition of $\ell$-smoothness in \citet{li2023generalized,tyurin2025unified}. For $\ell\equiv L$ we obtain the class of convex functions with $L$-Lipschitz gradient. The class $\F_\ell^d$ is invariant under orthogonal changes of coordinates, shifts of the argument and addition of a constant.

Checking condition~\eqref{eq:class} directly is inconvenient. The following lemma shows that it suffices to check it along lines, and that for twice continuously differentiable functions it is equivalent to the familiar Hessian bound. Its proof is given in Appendix~\ref{app:tools-class}.

\proofdown{lem:class}
\begin{lemma}[Characterization of the class]\label{lem:class}
Let $f:\R^d\to\R$ be convex and differentiable.
\begin{enumerate}
\item[(a)] $f\in\F_\ell^d$ if and only if for all $x,v\in\R^d$ the map $t\mapsto\nabla f(x+tv)$ is locally absolutely continuous on $\R$ and for almost all $t$
\begin{equation}\label{eq:line-derivative}
 \Bigl\|\frac{d}{dt}\nabla f(x+tv)\Bigr\|\le\ell\bigl(\norm{\nabla f(x+tv)}\bigr)\norm v .
\end{equation}
\item[(b)] If $f\in C^2(\R^d)$, then $f\in\F_\ell^d$ if and only if for all $x$
\begin{equation}\label{eq:hessian-bound}
 0\preceq\nabla^2f(x)\preceq\ell\bigl(\norm{\nabla f(x)}\bigr)I .
\end{equation}
\end{enumerate}
\end{lemma}

By shifting the argument we may assume $x^0=0$, which we do. Two problem classes are considered.
\begin{itemize}
\item $\mathcal P_{\rm R}$ consists of all functions $f\in\F_{L_0,L_1}^d$, $d\ge1$, that have a minimizer $x^\ast$ with $\norm{x^\ast}\le R$ (\emph{the radius-only class}).
\item $\mathcal P_\Delta$ consists of the functions from $\mathcal P_{\rm R}$ that additionally satisfy $f(0)-f^\ast\le\Delta$ (\emph{the class with bounded initial optimality gap}).
\end{itemize}
A method \emph{solves} a class with $N$ oracle calls in dimension $d$ if its output satisfies $f(\hat x)-f^\ast\le\varepsilon$ for all functions of the class on $\R^d$. The \emph{oracle complexity} $\mathcal N(\mathcal P)$ of a class $\mathcal P$ is the smallest $N$ such that for every dimension $d$ there exists a method solving all problems of the class in this dimension with at most $N$ calls. The precise model of a deterministic method is needed for the lower bound and is given in Section~\ref{sec:lower}. We are primarily interested in the complexity $\mathcal N(\mathcal P_{\rm R})$, and the class $\mathcal P_\Delta$ is considered in Appendix~\ref{app:gap}.

\section{Upper bound}\label{sec:upper}

In this section we construct a method that achieves the upper bound on the complexity of the class $\mathcal P_{\rm R}$. The method consists of two stages.
\begin{itemize}
\item \emph{The first stage} (the method SHERPA) passes through the large-gradient regime. Starting from $x^0$, it finds a point $x'$ such that
\begin{equation}\label{eq:stage-one-goal}
\begin{aligned}
 f(x')-f^\ast&\le\frac{L_0}{128L_1^2},\\
 \norm{x'-x^\ast}&\le R .
\end{aligned}
\end{equation}
This is an accelerated method based on Nesterov's framework of estimating sequences \citep[Section~2.2]{nesterov2018lectures}, with an inexact line search.
\item \emph{The second stage} (the method CRUISE) is started from $x'$. In a neighborhood of such a point, a function of the class behaves like a function with $2L_0$-Lipschitz gradient, and we run the standard accelerated method with a constant stepsize $1/(2L_0)$, using the analysis of \citet{tyurin2026accelerated}.
\end{itemize}
Subsections~\ref{sec:upper-descent}--\ref{sec:count3} are devoted to the first stage, Subsection~\ref{sec:upper-second} to the second one, and Subsection~\ref{sec:upper-proof} proves the final bound. Throughout this section $f\in\F_{L_0,L_1}^d$ and $x^\ast$ is a minimizer of $f$. The proofs of the auxiliary lemmas are deferred to Appendix~\ref{app:tools}.

\subsection{Descent lemma}\label{sec:upper-descent}

The first lemma is an analogue of the classical descent lemma \citep[Lemma~1.2.3]{nesterov2018lectures}. The gradient step is taken with a length computed from the local smoothness constant $L_0+L_1\norm{\nabla f(y)}$.

\proofdown{lem:safe-step}
\begin{lemma}[Descent lemma]\label{lem:safe-step}
For every point $y\in\R^d$, the point
\[
 y^+=y-\frac{\nabla f(y)}{4\bigl(L_0+L_1\norm{\nabla f(y)}\bigr)}
\]
satisfies
\begin{gather*}
 f(y^+)\le f(y)-\frac{5\norm{\nabla f(y)}^2}{24\bigl(L_0+L_1\norm{\nabla f(y)}\bigr)},\\
 \norm{y^+-x^\ast}^2\le\norm{y-x^\ast}^2-\frac{\norm{\nabla f(y)}^2}{24\bigl(L_0+L_1\norm{\nabla f(y)}\bigr)^2} .
\end{gather*}
\end{lemma}

The proof of Lemma~\ref{lem:safe-step} is given in Appendix~\ref{app:tools-descent}. The length of such a step is at most $1/(4L_1)$. Over this distance, the local smoothness constant changes by a factor of at most $e^{1/4}$ (Lemma~\ref{lem:gronwall}), so the function behaves like a function with Lipschitz gradient.

\subsection{Inequality at the minimizer}\label{sec:upper-anchored}

For $G\ge0$ we introduce the function
\begin{equation}\label{eq:Phi}
 \Phi(G)=\frac{L_0}{L_1^2}\Bigl[\Bigl(1+\frac{L_1G}{L_0}\Bigr)\log\Bigl(1+\frac{L_1G}{L_0}\Bigr)-\frac{L_1G}{L_0}\Bigr] .
\end{equation}
It is convex, $\Phi(0)=0$ and
\begin{equation}\label{eq:Phi-derivatives}
\begin{aligned}
 &\Phi'(G)=\frac1{L_1}\log\Bigl(1+\frac{L_1G}{L_0}\Bigr),\\
 &\Phi''(G)=\frac1{L_0+L_1G} .
\end{aligned}
\end{equation}
For small $G$ the function $\Phi(G)$ behaves like $G^2/(2L_0)$, and for large $G$ like $(G/L_1)\log(L_1G/L_0)$.

By convexity, the affine approximation of $f$ at $y$ underestimates $f^\ast$ at the point $x^\ast$. In the class $\F_{L_0,L_1}^d$ it underestimates $f^\ast$ by at least $\Phi(\norm{\nabla f(y)})$.

\proofdown{lem:anchored}
\begin{lemma}[Inequality at the minimizer]\label{lem:anchored}
For every point $y\in\R^d$
\[
 f(y)+\ip{\nabla f(y)}{x^\ast-y}+\Phi\bigl(\norm{\nabla f(y)}\bigr)\le f^\ast .
\]
\end{lemma}

The proof of Lemma~\ref{lem:anchored} is given in Appendix~\ref{app:tools-anchored}. For twice continuously differentiable functions, this is a special case of the inequality of \citet[Lemma~2.7]{vankov2025optimizing}. Thanks to the term $\Phi(\norm{\nabla f(y)})$, the first-stage method does not need a bound on the initial optimality gap.

\subsection{Inexact line search}\label{sec:upper-search}

At every iteration, the accelerated method needs a point $y$ on the segment $[v,x]$ at which simultaneously $f(y)\le f(x)$ and the quantity $\ip{\nabla f(y)}{y-v}$ is small. Exact line search would give $\ip{\nabla f(y)}{y-v}\le0$, but it is not guaranteed in a finite number of oracle calls. Bisection with a number of steps fixed in advance suffices (Algorithm~\ref{alg:search}). We call this method ICEPICK (Inexact Convex Edge Probing by Interval Cutting with Known Query Bound).

\begin{algorithm}[ht]
\caption{ICEPICK: inexact line search}\label{alg:search}
\begin{algorithmic}[1]
\Require points $v,x\in\R^d$, numbers $r,\delta>0$, a positive integer $m$
\State Query $\nabla f(v)$ and $\nabla f(x)$
\If{$r\norm{\nabla f(v)}\le\delta$} \Comment{stopping criterion}
 \State \Return $y=v$
\ElsIf{$r\norm{\nabla f(x)}\le\delta$} \Comment{stopping criterion}
 \State \Return $y=x$
\ElsIf{$\ip{\nabla f(v)}{x-v}\ge0$}
 \State \Return $y=v$
\ElsIf{$\ip{\nabla f(x)}{x-v}\le0$}
 \State \Return $y=x$
\EndIf
\State $t_-=0$, $t_+=1$
\For{$j=1,\dots,m$}
 \State $t=(t_-+t_+)/2$, query $\nabla f\bigl(v+t(x-v)\bigr)$
 \If{$r\norm{\nabla f(v+t(x-v))}\le\delta$} \Comment{stopping criterion}
  \State \Return $y=v+t(x-v)$
 \ElsIf{$\ip{\nabla f(v+t(x-v))}{x-v}=0$}
  \State \Return $y=v+t(x-v)$
 \ElsIf{$\ip{\nabla f(v+t(x-v))}{x-v}<0$}
  \State $t_-=t$
 \Else
  \State $t_+=t$
 \EndIf
\EndFor
\State \Return $y=v+t_+(x-v)$
\end{algorithmic}
\end{algorithm}

\proofdown{lem:search}
\begin{lemma}[Inexact line search]\label{lem:search}
Let $r,\delta>0$, let the points $v$ and $x$ lie in the ball $B(x^\ast,r)$, and let a positive integer $m$ satisfy
\begin{equation}\label{eq:search-m}
 2^m\ge128\Bigl(1+\frac{L_0r^2}\delta+L_1r\Bigr)^2 .
\end{equation}
Then Algorithm~\ref{alg:search} makes at most $m+2$ oracle calls at points of the segment $[v,x]$ and returns a point $y\in[v,x]$ for which one of the following two conditions holds.
\begin{enumerate}
\item[(a)] The stopping criterion was triggered. Then $f(y)-f^\ast\le\delta$.
\item[(b)] The inequalities
\begin{equation}\label{eq:search}
\begin{aligned}
 f(y)&\le f(x),\\
 \ip{\nabla f(y)}{y-v}&\le\frac{\norm{\nabla f(y)}^2}{16\bigl(L_0+L_1\norm{\nabla f(y)}\bigr)}
\end{aligned}
\end{equation}
hold.
\end{enumerate}
\end{lemma}

The proof of Lemma~\ref{lem:search} is given in Appendix~\ref{app:tools-search}. The number of steps $m$ depends logarithmically on $L_0r^2/\delta$ and $L_1r$. This is exactly where the logarithmic factor in the upper bound comes from. Algorithm~\ref{alg:search} uses only gradients: it does not need function values.

\subsection{First-stage method}\label{sec:upper-method}

The first-stage method SHERPA (Smoothness-aware High-gradient Escape with Radius-Preserving Acceleration, Algorithm~\ref{alg:am}) is designed as an accelerated method based on estimating sequences \citep[Section~2.2]{nesterov2018lectures} and is close to the AGMsDR method \citep[Algorithm~1]{vankov2025optimizing}. It maintains two sequences of points: $x^k$ with decreasing function values and $v^k$, from which a lower model of the function is built. Here and below, the sum of the first $k$ terms of a sequence is denoted by the corresponding capital letter: $A_k=a_1+\dots+a_k$, $A_0=0$. The coefficient $a_k>0$ is defined by the equation
\begin{equation}\label{eq:coef}
 a_k^2=\frac{A_{k-1}+a_k}{7\bigl(L_0+L_1\norm{\nabla f(y^k)}\bigr)}+\frac{a_k\,\Phi\bigl(\norm{\nabla f(y^k)}\bigr)}{\norm{\nabla f(y^k)}^2} .
\end{equation}
This is a quadratic equation in $a_k$ with a unique positive root. Without the last term, it coincides with the usual coefficient rule of the accelerated method in which the Lipschitz constant of the gradient is replaced by $7(L_0+L_1\norm{\nabla f(y^k)})$. The last term corresponds to the additional term $\Phi$ of Lemma~\ref{lem:anchored}. It increases the coefficient at iterations with a large gradient.

\begin{algorithm}[ht]
\caption{SHERPA (first stage): accelerated method with inexact line search}\label{alg:am}
\begin{algorithmic}[1]
\Require starting point $x^0$, parameters $L_0$, $L_1$, $R$ with $L_1R\ge1$
\State $m=\Bigl\lceil\log_2\bigl(128(1+128L_1^2R^2+L_1R)^2\bigr)\Bigr\rceil$, $K=\lceil60L_1R\rceil$
\State $v^0=x^0$, $A_0=0$
\For{$k=1,\dots,K$}
 \State $y^k$ is the output of Algorithm~\ref{alg:search} on the segment $[v^{k-1},x^{k-1}]$ with $r=R$, $\delta=L_0/(128L_1^2)$ and $m$
 \If{the stopping criterion was triggered}
  \State \Return $y^k$
 \EndIf
 \State $a_k$ is the positive root of equation~\eqref{eq:coef}, $A_k=A_{k-1}+a_k$
 \State $x^k=y^k-\dfrac{\nabla f(y^k)}{4(L_0+L_1\norm{\nabla f(y^k)})}$
 \State $v^k=v^{k-1}-a_k\nabla f(y^k)$
\EndFor
\State \Return $x^K$
\end{algorithmic}
\end{algorithm}

Equation~\eqref{eq:coef} and the inequality $A_{k-1}\ge0$ imply two properties of the coefficients:
\begin{equation}\label{eq:coef-basic}
\begin{aligned}
 a_k&\le A_k,\\
 a_k^2&\ge\frac{A_k}{7\bigl(L_0+L_1\norm{\nabla f(y^k)}\bigr)} .
\end{aligned}
\end{equation}
At every iteration, Algorithm~\ref{alg:search} makes at most $m+2$ oracle calls, and the points $x^k$ and $v^k$ are computed without oracle calls. Therefore, $K$ iterations cost at most $K(m+2)$ calls. The method does not need function values.

\subsection{Optimality gap and distance to the minimizer}\label{sec:upper-energy}

The following lemma simultaneously bounds the optimality gap, the weighted sum of squared gradient norms and the distance from $v^k$ to the minimizer.

\begin{lemma}\label{lem:energy}
Let $\norm{x^0-x^\ast}\le R$, and suppose that Algorithm~\ref{alg:am} completes iterations $1,\dots,k$ without meeting the stopping criterion. Then
\begin{equation}\label{eq:energy}
 A_k\bigl(f(x^k)-f^\ast\bigr)+\frac12\sum_{i=1}^ka_i^2\norm{\nabla f(y^i)}^2+\frac12\norm{v^k-x^\ast}^2\le\frac{R^2}2 ,
\end{equation}
and all points $x^i$, $v^i$, $y^i$, $i\le k$, lie in the ball $B(x^\ast,R)$.
\end{lemma}

\begin{proof}
Induction on $k$. For $k=0$, inequality~\eqref{eq:energy} coincides with the condition $\norm{x^0-x^\ast}\le R$. Suppose the statement holds for $k-1$.
\begin{itemize}
\item \emph{The iteration is well defined.} By the induction hypothesis, the points $v^{k-1}$ and $x^{k-1}$ lie in the ball $B(x^\ast,R)$. The ball is convex, so the segment $[v^{k-1},x^{k-1}]$ and the point $y^k$ lie in it, and Lemma~\ref{lem:search} applies. The stopping criterion was not triggered, so inequalities~\eqref{eq:search} hold. Together with Lemma~\ref{lem:safe-step} they give $x^k\in B(x^\ast,R)$ and
\begin{equation}\label{eq:energy-step}
\begin{aligned}
 &f(y^k)\le f(x^{k-1}),\\
 &\ip{\nabla f(y^k)}{y^k-v^{k-1}}\le\frac{\norm{\nabla f(y^k)}^2}{16\bigl(L_0+L_1\norm{\nabla f(y^k)}\bigr)},\\
 &f(x^k)\le f(y^k)-\frac{5\norm{\nabla f(y^k)}^2}{24\bigl(L_0+L_1\norm{\nabla f(y^k)}\bigr)} .
\end{aligned}
\end{equation}
\item \emph{Distance.} Since $v^k=v^{k-1}-a_k\nabla f(y^k)$,
\begin{equation}\label{eq:energy-dist}
\begin{aligned}
 \MoveEqLeft\frac12\norm{v^k-x^\ast}^2-\frac12\norm{v^{k-1}-x^\ast}^2\\
 &=a_k\ip{\nabla f(y^k)}{y^k-v^{k-1}}-a_k\ip{\nabla f(y^k)}{y^k-x^\ast}+\frac12a_k^2\norm{\nabla f(y^k)}^2 .
\end{aligned}
\end{equation}
We bound the first two terms on the right:
\begin{align*}
 a_k\ip{\nabla f(y^k)}{y^k-v^{k-1}}
 &\bypad{\le}{\eqref{eq:energy-step}}{\text{Lemma~\ref{lem:anchored}}}\frac{a_k\norm{\nabla f(y^k)}^2}{16\bigl(L_0+L_1\norm{\nabla f(y^k)}\bigr)},\\
 -a_k\ip{\nabla f(y^k)}{y^k-x^\ast}
 &\by{\le}{\text{Lemma~\ref{lem:anchored}}}-a_k\bigl(f(y^k)-f^\ast\bigr)-a_k\Phi\bigl(\norm{\nabla f(y^k)}\bigr).
\end{align*}
\item \emph{Optimality gap.} Since $A_k=A_{k-1}+a_k$,
\begin{align*}
 \MoveEqLeft A_k\bigl(f(x^k)-f^\ast\bigr)-A_{k-1}\bigl(f(x^{k-1})-f^\ast\bigr)\\
 &\by{\le}{\eqref{eq:energy-step}}A_k\bigl(f(y^k)-f^\ast\bigr)-\frac{5A_k\norm{\nabla f(y^k)}^2}{24\bigl(L_0+L_1\norm{\nabla f(y^k)}\bigr)}-A_{k-1}\bigl(f(y^k)-f^\ast\bigr)\\
 &\bypad{=}{}{\eqref{eq:energy-step}}a_k\bigl(f(y^k)-f^\ast\bigr)-\frac{5A_k\norm{\nabla f(y^k)}^2}{24\bigl(L_0+L_1\norm{\nabla f(y^k)}\bigr)} .
\end{align*}
\item \emph{Sum.} Multiplying equation~\eqref{eq:coef} by $\norm{\nabla f(y^k)}^2$, we obtain
\[
 a_k\Phi\bigl(\norm{\nabla f(y^k)}\bigr)=a_k^2\norm{\nabla f(y^k)}^2-\frac{A_k\norm{\nabla f(y^k)}^2}{7\bigl(L_0+L_1\norm{\nabla f(y^k)}\bigr)} .
\]
We add~\eqref{eq:energy-dist}, the two bounds on its terms and the bound on the optimality gap, and then substitute this equality:
\begin{align*}
 \MoveEqLeft\Bigl[A_k\bigl(f(x^k)-f^\ast\bigr)+\frac12\norm{v^k-x^\ast}^2\Bigr]-\Bigl[A_{k-1}\bigl(f(x^{k-1})-f^\ast\bigr)+\frac12\norm{v^{k-1}-x^\ast}^2\Bigr]\\
 &\bypad{\le}{}{\eqref{eq:coef-basic}}\Bigl(\frac{a_k}{16}-\frac{5A_k}{24}\Bigr)\frac{\norm{\nabla f(y^k)}^2}{L_0+L_1\norm{\nabla f(y^k)}}-a_k\Phi\bigl(\norm{\nabla f(y^k)}\bigr)+\frac12a_k^2\norm{\nabla f(y^k)}^2\\
 &\bypad{=}{}{\eqref{eq:coef-basic}}\Bigl(\frac{a_k}{16}-\frac{5A_k}{24}+\frac{A_k}7\Bigr)\frac{\norm{\nabla f(y^k)}^2}{L_0+L_1\norm{\nabla f(y^k)}}-\frac12a_k^2\norm{\nabla f(y^k)}^2\\
 &\by{\le}{\eqref{eq:coef-basic}}-\frac{A_k\norm{\nabla f(y^k)}^2}{336\bigl(L_0+L_1\norm{\nabla f(y^k)}\bigr)}-\frac12a_k^2\norm{\nabla f(y^k)}^2\\
 &\bypad{\le}{}{\eqref{eq:coef-basic}}-\frac12a_k^2\norm{\nabla f(y^k)}^2 .
\end{align*}
Adding this inequality to inequality~\eqref{eq:energy} for $k-1$, we obtain~\eqref{eq:energy} for $k$. It implies $\norm{v^k-x^\ast}\le R$, which completes the induction.\qedhere
\end{itemize}
\end{proof}

\subsection{Number of first-stage iterations}\label{sec:count3}

The initial optimality gap is unknown and may be exponentially large in $L_1R$, so a bound in terms of this gap does not give the desired guarantee in terms of the radius. Instead, we use the two quantities bounded by Lemma~\ref{lem:energy}: the total weight $A_K$ and the sum of the squares $a_i^2\norm{\nabla f(y^i)}^2$. Iterations with a very large gradient increase the second sum, so there are few of them. At the remaining iterations the weight $A_i$ is not too small, and the number of such iterations is bounded in the same way as for the accelerated method on functions with Lipschitz gradient.

\begin{lemma}[Number of first-stage iterations]\label{lem:count-anchor}
Let $\norm{x^0-x^\ast}\le R$ and $L_1R\ge1$. If Algorithm~\ref{alg:am} completes $K$ iterations without meeting the stopping criterion and
\[
 f(x^K)-f^\ast>\frac{L_0}{128L_1^2},
\]
then $K<60L_1R$.
\end{lemma}

\begin{proof}
\leavevmode
\begin{itemize}
\item \emph{Two sums.} By Lemma~\ref{lem:energy} and the condition on the optimality gap
\begin{equation}\label{eq:count-anchor-budget}
\begin{aligned}
 &L_0A_K<64L_1^2R^2,\\
 &\sum_{i=1}^Ka_i^2\norm{\nabla f(y^i)}^2\le R^2 .
\end{aligned}
\end{equation}
\item \emph{Bounds for one iteration.} From inequalities~\eqref{eq:coef-basic}
\begin{equation}\label{eq:count-step}
\begin{gathered}
\begin{aligned}
 \sqrt{L_0A_i}-\sqrt{L_0A_{i-1}}=\frac{\sqrt{L_0}\,a_i}{\sqrt{A_i}+\sqrt{A_{i-1}}}
 &\by{\ge}{A_{i-1}\le A_i}\frac{\sqrt{L_0}\,a_i}{2\sqrt{A_i}}\\
 &\bypad{\ge}{\eqref{eq:coef-basic}}{A_{i-1}\le A_i}\frac1{2\sqrt7}\sqrt{\frac{L_0}{L_0+L_1\norm{\nabla f(y^i)}}},
\end{aligned}\\
\begin{aligned}
 L_0A_i=\frac{L_0A_i^2}{A_i}
 &\by{\ge}{\eqref{eq:coef-basic}}\frac{L_0a_i^2}{A_i}\\
 &\by{\ge}{\eqref{eq:coef-basic}}\frac{L_0}{7\bigl(L_0+L_1\norm{\nabla f(y^i)}\bigr)} .
\end{aligned}
\end{gathered}
\end{equation}
Divide equation~\eqref{eq:coef} by $a_i$ and multiply by $L_1\norm{\nabla f(y^i)}$:
\begin{equation}\label{eq:anchor-identity}
 L_1a_i\norm{\nabla f(y^i)}
 =\frac{L_1\Phi(\norm{\nabla f(y^i)})}{\norm{\nabla f(y^i)}}
 +\frac{A_i}{7a_i}\cdot\frac{L_1\norm{\nabla f(y^i)}}{L_0+L_1\norm{\nabla f(y^i)}} .
\end{equation}
Both terms on the right are nonnegative. By definition~\eqref{eq:Phi}
\begin{equation}\label{eq:Phi-ratio}
\begin{aligned}
 \frac{L_1\Phi(\norm{\nabla f(y^i)})}{\norm{\nabla f(y^i)}}
 &=\Bigl(1+\frac{L_0}{L_1\norm{\nabla f(y^i)}}\Bigr)\log\Bigl(1+\frac{L_1\norm{\nabla f(y^i)}}{L_0}\Bigr)-1\\
 &\ge\log\Bigl(1+\frac{L_1\norm{\nabla f(y^i)}}{L_0}\Bigr)-1 .
\end{aligned}
\end{equation}
Multiplying~\eqref{eq:anchor-identity} by $7a_i/A_i$ and adding the first inequality of~\eqref{eq:count-step}, we obtain
\begin{equation}\label{eq:count-two-groups}
\begin{aligned}
 \MoveEqLeft2\sqrt7\Bigl(\sqrt{L_0A_i}-\sqrt{L_0A_{i-1}}\Bigr)+\Bigl(\frac{7L_1a_i^2\norm{\nabla f(y^i)}}{A_i}\Bigr)^{2/3}\\
 &\ge\Bigl(\frac{L_0}{L_0+L_1\norm{\nabla f(y^i)}}\Bigr)^{1/2}+\Bigl(\frac{L_1\norm{\nabla f(y^i)}}{L_0+L_1\norm{\nabla f(y^i)}}\Bigr)^{2/3}\\
 &\by{\ge}{(\mathrm a)}\frac{L_0}{L_0+L_1\norm{\nabla f(y^i)}}+\frac{L_1\norm{\nabla f(y^i)}}{L_0+L_1\norm{\nabla f(y^i)}}=1 .
\end{aligned}
\end{equation}
Here step (a) uses the inequalities $t^{1/2}\ge t$ and $t^{2/3}\ge t$ for $0\le t\le1$.
\item \emph{Iterations with a large product $L_1a_i\norm{\nabla f(y^i)}$.} Let $J$ be the set of indices $i\le K$ with $L_1a_i\norm{\nabla f(y^i)}\ge\sqrt{L_1R}$. By the second inequality of~\eqref{eq:count-anchor-budget}
\[
 |J|\,L_1R\le\sum_{i\in J}L_1^2a_i^2\norm{\nabla f(y^i)}^2\le L_1^2R^2 ,
\]
that is, $|J|\le L_1R$.
\item \emph{Remaining iterations.} Let $I$ be the set of the remaining indices. If it is empty, then $K=|J|\le L_1R$. Otherwise, for $i\in I$, by~\eqref{eq:anchor-identity} and~\eqref{eq:Phi-ratio}
\[
 \log\Bigl(1+\frac{L_1\norm{\nabla f(y^i)}}{L_0}\Bigr)-1\le L_1a_i\norm{\nabla f(y^i)}<\sqrt{L_1R} ,
\]
and the second inequality of~\eqref{eq:count-step} gives
\[
 L_0A_i\ge\frac17\Bigl(1+\frac{L_1\norm{\nabla f(y^i)}}{L_0}\Bigr)^{-1}>\frac17e^{-\sqrt{L_1R}-1} .
\]
Let $j$ be the smallest index in $I$. The weights $A_i$ increase, so for $i>j$ we have $A_{i-1}>0$ and
\[
 \frac{a_i}{A_i}\le-\log\Bigl(1-\frac{a_i}{A_i}\Bigr)=\log\frac{A_i}{A_{i-1}} .
\]
Together with the inequality $a_j\le A_j$ this gives
\begin{align*}
 \sum_{i\in I}\frac{a_i}{A_i}
 &\le1+\log\frac{A_K}{A_j}
 \by{<}{\eqref{eq:count-anchor-budget}}1+\log\bigl(448L_1^2R^2e^{\sqrt{L_1R}+1}\bigr)\\
 &\bypad{=}{}{\eqref{eq:count-anchor-budget}}\sqrt{L_1R}+2+\log448+2\log(L_1R)\\
 &\bypad{<}{(\mathrm b)}{\eqref{eq:count-anchor-budget}}5\sqrt{L_1R}+\frac92
 \by{\le}{L_1R\ge1}\frac{19}2\sqrt{L_1R} .
\end{align*}
Here step (b) uses the inequalities $\log448<13/2$ and $\log t\le2\sqrt t-2$ for $t\ge1$.
Sum inequalities~\eqref{eq:count-two-groups} over $i\in I$. The sum of the first terms is at most $2\sqrt7\sqrt{L_0A_K}<16\sqrt7L_1R$. To the second terms we apply Hölder's inequality with exponents $3$ and $3/2$:
\begin{align*}
 |I|
 &<16\sqrt7L_1R+7^{2/3}\sum_{i\in I}\bigl(L_1a_i\norm{\nabla f(y^i)}\bigr)^{2/3}\Bigl(\frac{a_i}{A_i}\Bigr)^{2/3}\\
 &\bypad{\le}{}{\eqref{eq:count-anchor-budget}}16\sqrt7L_1R+7^{2/3}\Bigl(\sum_{i\in I}L_1^2a_i^2\norm{\nabla f(y^i)}^2\Bigr)^{1/3}\Bigl(\sum_{i\in I}\frac{a_i}{A_i}\Bigr)^{2/3}\\
 &\by{<}{\eqref{eq:count-anchor-budget}}16\sqrt7L_1R+7^{2/3}(L_1R)^{2/3}\Bigl(\frac{19}2\sqrt{L_1R}\Bigr)^{2/3}\\
 &\bypad{=}{}{\eqref{eq:count-anchor-budget}}\Bigl(16\sqrt7+\Bigl(\frac{133}2\Bigr)^{2/3}\Bigr)L_1R\\
 &\bypad{<}{(\mathrm c)}{\eqref{eq:count-anchor-budget}}59L_1R .
\end{align*}
Here step (c) uses the inequalities $16\sqrt7<42.4$ and $(133/2)^{2/3}<16.5$, which are equivalent to $1792<42.4^2$ and $133^2<4\cdot16.5^3$.
\item \emph{Conclusion.} Every iteration belongs to exactly one of the sets $J$ and $I$, therefore
\[
 K=|J|+|I|<L_1R+59L_1R=60L_1R .\qedhere
\]
\end{itemize}
\end{proof}

\subsection{Second stage}\label{sec:upper-second}

At the second stage, the accelerated gradient method with a constant stepsize $1/(2L_0)$ \citep[Algorithm~1]{tyurin2026accelerated} is run (Algorithm~\ref{alg:agm}). Its structure is close to the similar triangles method \citep{gasnikov2018universal}. We call the variant that we use CRUISE (Constant-step Refinement Using Inertia in Smooth Environments). This is the name of a variant of the accelerated method of \citet{tyurin2026accelerated} started from the point found at the first stage, not a new method developed by the authors of this paper. Each iteration computes one new gradient.

\begin{algorithm}[ht]
\caption{CRUISE (second stage): accelerated method with a constant stepsize from \citet{tyurin2026accelerated}}\label{alg:agm}
\begin{algorithmic}[1]
\Require starting point $x'$, parameters $L_0$, $L_1$, $R$, number of iterations $N$
\State $y^0=w^0=x'$, $\Gamma_0=L_0/(64L_1^2R^2)$
\For{$k=0,\dots,N-1$}
 \State $\lambda_k=\sqrt{\Gamma_k/(2L_0)}$
 \State $y^{k+1}=\dfrac1{1+\lambda_k}\Bigl(y^k+\lambda_kw^k-\dfrac{\nabla f(y^k)}{2L_0}\Bigr)$
 \State $w^{k+1}=w^k-\dfrac{\lambda_k}{\Gamma_k}\nabla f(y^{k+1})$
 \State $\Gamma_{k+1}=\Gamma_k/(1+\lambda_k)$
\EndFor
\State \Return $y^N$
\end{algorithmic}
\end{algorithm}

\begin{proposition}[Second stage]\label{prop:agm}
Let $L_1R\ge1$, and let the point $x'$ satisfy
\begin{align*}
 f(x')-f^\ast&\le\frac{L_0}{128L_1^2},\\
 \norm{x'-x^\ast}&\le R .
\end{align*}
Then the points of Algorithm~\ref{alg:agm} satisfy, for all $k\ge1$,
\[
 f(y^k)-f^\ast\le\frac{18L_0R^2}{k^2} .
\]
In particular, for
\begin{equation}\label{eq:agm-N}
 N=\Bigl\lceil\sqrt{\frac{18L_0R^2}\varepsilon}\,\Bigr\rceil
\end{equation}
the algorithm finds a point with optimality gap at most $\varepsilon$ after at most
\[
 5\sqrt{\frac{L_0R^2}\varepsilon}+2
\]
oracle calls.
\end{proposition}

\begin{proof}
This is Theorem~3.2 of \citet{tyurin2026accelerated} for $\ell(s)=L_0+L_1s$ with the accuracy parameter equal to twice the right-hand side of the first inequality of~\eqref{eq:stage-one-goal}. We verify its assumptions.
\begin{itemize}
\item The author of \citet{tyurin2026accelerated} states the result for twice differentiable functions, but notes (ibid., after Lemma~2.2) that instead of the Hessian bound the inequality
\[
 \norm{\nabla f(z)-\nabla f(y)}\le\bigl(L_0+L_1\norm{\nabla f(y)}\bigr)\frac{e^{L_1\norm{z-y}}-1}{L_1}
\]
suffices. For our class this is Lemma~\ref{lem:gronwall}(b).
\item The function in the assumptions of Theorem~3.2 is strictly increasing on $[0,\infty)$, since for $t>0$
\[
 \frac{d}{dt}\,\frac{t^2}{2(L_0+4L_1t)}=\frac{t(L_0+2L_1t)}{(L_0+4L_1t)^2}>0 .
\]
\item The condition of Theorem~3.2 on the accuracy parameter holds with equality:
\[
 \ell\Bigl(8\sqrt{\ell(0)\,\frac{L_0}{64L_1^2}}\Bigr)=L_0+8L_1\sqrt{\frac{L_0^2}{64L_1^2}}=2L_0=2\ell(0) .
\]
\item In Theorem~3.2 the starting point of the accelerated method is obtained by gradient descent, but the proof (ibid., Appendix~E.1, formula~(32)) uses only its two properties~\eqref{eq:stage-one-goal}. For the point $x'$ they hold by assumption.
\item The number of preliminary iterations in Theorem~3.2 is zero:
\[
 \frac{\Gamma_0}{8L_0}=\frac1{512L_1^2R^2}\le\frac49 ,
\]
therefore
\[
 \max\Bigl\{1+\frac12\log_{3/2}\frac{\Gamma_0}{8L_0},\ 0\Bigr\}=0 .
\]
\end{itemize}
Theorem~3.2 gives the bound in the statement of the proposition. With $N$ chosen as in~\eqref{eq:agm-N}, its right-hand side for $k=N$ is at most $\varepsilon$. The algorithm computes gradients at the points $y^0,\dots,y^N$, so the number $N+1$ of its oracle calls satisfies
\[
 N+1\le\sqrt{\frac{18L_0R^2}\varepsilon}+2\le5\sqrt{\frac{L_0R^2}\varepsilon}+2 .
\]
\end{proof}

\subsection{Final bound}\label{sec:upper-proof}

\begin{theorem}[Upper bound for the radius-only class]\label{thm:upper}
For all $L_0,L_1,R,\varepsilon>0$
\[
 \mathcal N(\mathcal P_{\rm R})\le5\sqrt{\frac{L_0R^2}\varepsilon}+240\,(1+L_1R)\log_2\bigl(64(1+L_1R)\bigr)+2 .
\]
The method achieving this bound uses only gradients and the numbers $L_0$, $L_1$, $R$, $\varepsilon$. It needs neither the value $f^\ast$, nor a bound on the initial optimality gap, nor function values.
\end{theorem}

\begin{proof}
\leavevmode
\begin{itemize}
\item \emph{Reduction to the case $L_1R\ge1$.} If $L_1R<1$, then the class $\F_{L_0,L_1}^d$ is contained in the class $\F_{L_0,1/R}^d$. We then replace $L_1$ by $1/R$ in all parameters of both stages. After the replacement $K=60$ and $m=22$, so the first stage makes at most $60\cdot24=240\log_264$ oracle calls, which does not exceed $240(1+L_1R)\log_2\bigl(64(1+L_1R)\bigr)$. From now on $L_1R\ge1$.
\item \emph{First stage.} Run Algorithm~\ref{alg:am}. Its number $m$ satisfies condition~\eqref{eq:search-m} with $r=R$ and $\delta=L_0/(128L_1^2)$. Let $x'$ be the returned point. We prove that it satisfies inequalities~\eqref{eq:stage-one-goal}. If the stopping criterion was triggered at iteration $k$, then the first of them is given by Lemma~\ref{lem:search}. Moreover, $x'$ lies on the segment $[v^{k-1},x^{k-1}]$, whose endpoints lie in the ball $B(x^\ast,R)$ by Lemma~\ref{lem:energy}. Otherwise all $K\ge60L_1R$ iterations are performed, the first inequality is given by Lemma~\ref{lem:count-anchor}, and the second one by Lemma~\ref{lem:energy}.
\item \emph{Number of first-stage calls.} It is at most $K(m+2)$. For $L_1R\ge1$
\begin{align*}
 &K\le60L_1R+1\le60(1+L_1R),\\
 &m+2\by{\le}{(\mathrm a)}\log_2\bigl(2^{21}L_1^2R^2(1+L_1R)^2\bigr)+3\le4\log_2\bigl(64(1+L_1R)\bigr) .
\end{align*}
Here step (a) uses the inequality $1+128L_1^2R^2+L_1R\le128L_1R(1+L_1R)$. Hence, the first stage makes at most $240(1+L_1R)\log_2\bigl(64(1+L_1R)\bigr)$ oracle calls.
\item \emph{Second stage.} If $\varepsilon$ is not smaller than the right-hand side of the first inequality of~\eqref{eq:stage-one-goal}, we return $x'$. Otherwise we run Algorithm~\ref{alg:agm} from $x'$ with the number of iterations~\eqref{eq:agm-N}. By Proposition~\ref{prop:agm}, it finds a point with optimality gap at most $\varepsilon$, and its number of oracle calls is bounded as in this proposition.
\end{itemize}
Adding the numbers of calls at the two stages, we obtain the statement of the theorem. Both algorithms use only gradients.
\end{proof}

Thus,
\[
 \mathcal N(\mathcal P_{\rm R})=O\Bigl(\sqrt{\frac{L_0R^2}\varepsilon}+(1+L_1R)\log(2+L_1R)\Bigr).
\]
The classical term enters without a logarithm and with the constant~$5$, while the logarithmic factor affects only the $L_1R$ term and comes from the line search. Section~\ref{sec:lower} shows that this bound cannot be improved up to logarithmic factors.

\section{Lower bound}\label{sec:lower}

\paragraph{Oracle model.} We use the standard black-box model with a first-order oracle \citep{nemirovski1983problem}, \citep[Section~1.1.2]{nesterov2018lectures}. Given a query $z\in\R^d$, the oracle returns the pair $(f(z),\nabla f(z))$, and the method has no other information about the function. A \emph{deterministic method} making $N$ oracle calls in dimension $d$ is given by two rules: the $k$-th query point is an arbitrary function of the previous query points and the oracle answers at them, and the output $\hat x$ is an arbitrary function of all query points and answers. The method knows $L_0$, $L_1$, $R$, $\varepsilon$ (and the bound $\Delta$ on the initial optimality gap, if it is given), but does not know $f^\ast$ and $x^\ast$. Function values in the oracle answer only extend the capabilities of the method, so the lower bounds of this section hold a fortiori for methods that receive only gradients.

A lower bound means that for every deterministic method with fewer oracle calls there exists a function of the class whose optimality gap at the output exceeds $\varepsilon$. We use the classical scheme of \citet{nemirovski1983problem} in the form presented by \citet{carmon2020lower}. The scheme consists of two steps.
\begin{itemize}
\item \emph{Chain.} We construct a convex function on $\R^n$ for which every oracle call reveals at most one new coordinate. A Krylov method after $n-1$ calls does not leave the subspace where the last coordinate is zero, and there the optimality gap exceeds $\varepsilon$.
\item \emph{Orthogonal invariance principle.} The class is invariant under orthogonal transformations. Therefore, by rotating the chain in a space of dimension $2n$, the bound for Krylov methods carries over to arbitrary deterministic methods.
\end{itemize}
Both steps form a general tool, which is described in Subsection~\ref{sec:chains}. There we also construct chains conjugate to strongly convex functions on a polytope. Specific chains differ in the slopes of their links. In Subsection~\ref{sec:ordered} the slopes decrease geometrically, and the chain contains almost $L_1R$ links. In Appendix~\ref{app:gap} the slopes are equal, there are of order $(L_1R)^{2/3}$ links, but the initial optimality gap is small. The classical term independent of $L_1$ is given by the lower bound for functions with $L_0$-Lipschitz gradient \citep[Theorem~2.1.7]{nesterov2018lectures}, which we do not reprove. The main result of the section, Theorem~\ref{thm:lower-radius}, is proved in Subsection~\ref{sec:lower-radius-proof}.

\subsection{Chains and the orthogonal invariance principle}\label{sec:chains}

Let $e_1,\dots,e_n$ be the standard basis of $\R^n$. Set $\mathcal K_0=\{0\}$ and, for $j=1,\dots,n$,
\[
 \mathcal K_j=\operatorname{span}\{e_1,\dots,e_j\}.
\]

\begin{definition}[Chain]\label{def:zero-chain}
A differentiable function $\varphi:\R^n\to\R$ is called a \emph{chain} if for every $j=1,\dots,n$
\begin{equation}\label{eq:coord-chain}
 y\in\mathcal K_{j-1}\implies\nabla\varphi(y)\in\mathcal K_j .
\end{equation}
\end{definition}

These are first-order zero-chains in the terminology of \citet{carmon2020lower}, and the subspaces $\mathcal K_j$ play the role of Krylov subspaces. We call a method a \emph{Krylov method} if each of its query points and its output lie in the linear span of the gradients already obtained, that is, if it satisfies the linear span assumption \citep[Assumption~2.1.4]{nesterov2018lectures}. This extends the usual meaning of the term, in which Krylov subspaces are generated by a fixed matrix. After $n-1$ calls to a chain, such a method does not leave $\mathcal K_{n-1}$.

\paragraph{Chains on a polytope.} The quadratic chain does not use $L_1$. To obtain the term with $L_1R$, we need a chain whose curvature grows together with the gradient norm. We construct it as the conjugate of a strongly convex function $H$ on a polytope: the smaller the strong convexity modulus of $H$, the larger the curvature of the conjugate function. For a constant modulus, the gradient of the conjugate function is Lipschitz with constant equal to the reciprocal of the modulus \citep[Theorem~6.1.1]{nesterov2018lectures}. The following lemma allows the modulus to depend on the point.

\begin{lemma}[Smoothness of the conjugate function]\label{lem:conj}
Let $K\subset\R^n$ be a nonempty convex compact set, let the function $H$ be continuous on $K$, and suppose that for all $z,z'\in K$ its restriction to the segment $[z,z']$ is strongly convex with modulus
\[
 \frac1{L_0+L_1\max\{\norm z,\norm{z'}\}} .
\]
Then for every $\sigma\in\R^n$ the function
\[
 \varphi(y)=\max_{z\in K}\bigl\{\ip{y+\sigma}z-H(z)\bigr\}
\]
is convex and differentiable, its gradient $\nabla\varphi(y)$ is the unique maximizer $z(y)$, and for all $y,y'$
\begin{equation}\label{eq:conj-lip}
 \norm{\nabla\varphi(y)-\nabla\varphi(y')}\le\Bigl(L_0+L_1\max\bigl\{\norm{\nabla\varphi(y)},\norm{\nabla\varphi(y')}\bigr\}\Bigr)\norm{y-y'} .
\end{equation}
In particular, $\varphi\in\F_{L_0,L_1}^n$.
\end{lemma}

The modulus decreases as the norms of the endpoints of the segment grow. This is exactly why the curvature of $\varphi$ can grow together with the gradient norm. The gradient itself lies in $K$.

\begin{proof}
\leavevmode
\begin{itemize}
\item \emph{Gradient.} Denote $\Psi_y(z)=\ip{y+\sigma}z-H(z)$. The function $H$ is strictly convex on every segment, so $\Psi_y$ has a unique maximizer $z(y)$ on $K$. By Danskin's theorem, $\varphi$ is convex, differentiable and $\nabla\varphi(y)=z(y)$.
\item \emph{Growth around the maximizer.} Let $z=z(y)$, $z'\in K$, $z_t=z+t(z'-z)$, and let $m$ be the modulus from the assumption for the segment $[z,z']$. By assumption, the function $\Psi_y(z_t)$ remains concave in $t$ after adding the quadratic term with modulus $m$, therefore
\begin{align*}
 \Psi_y(z)
 &\by{\ge}{z=z(y)}\Psi_y(z_t)\\
 &\bypad{\ge}{(\mathrm a)}{z=z(y)}(1-t)\Psi_y(z)+t\Psi_y(z')+\frac m2t(1-t)\norm{z'-z}^2 ,
\end{align*}
where step (a) uses the strong concavity of $\Psi_y$ on the segment $[z,z']$ with modulus $m$.
Subtracting $(1-t)\Psi_y(z)$, dividing by $t$ and letting $t$ tend to zero, we obtain
\begin{equation}\label{eq:conj-growth}
 \Psi_y(z(y))-\Psi_y(z')\ge\frac m2\norm{z'-z(y)}^2 .
\end{equation}
\item \emph{Lipschitz continuity of the gradient.} Let $z=z(y)$, $z'=z(y')$, and let $m$ be the modulus for the segment $[z,z']$. Add~\eqref{eq:conj-growth} for the pair $(y,z')$ to the same inequality for the pair $(y',z)$. The values of $H$ cancel, and
\begin{align*}
 m\norm{z-z'}^2&\bypad{\le}{}{\text{C--S}}\ip{y-y'}{z-z'}\\
 &\by{\le}{\text{C--S}}\norm{y-y'}\,\norm{z-z'} .
\end{align*}
If $z=z'$, inequality~\eqref{eq:conj-lip} is obvious, otherwise dividing by $m\norm{z-z'}$ gives~\eqref{eq:conj-lip}. The maximum of the gradient norm over the endpoints of a segment does not exceed its supremum over the segment, so $\varphi\in\F_{L_0,L_1}^n$.\qedhere
\end{itemize}
\end{proof}

\paragraph{Orthogonal invariance principle.} The restriction on the method is removed by rotating the chain: the class is invariant under orthogonal transformations.

\begin{proposition}[Orthogonal invariance principle]\label{prop:reduction}
\leavevmode
\begin{enumerate}
\item[(a)] Let $\varphi:\R^n\to\R$ be a convex chain, $d\ge2n$, and let every $d\times n$ matrix $V$ with orthonormal columns $v_1,\dots,v_n$ be assigned the function
\[
 F_V(x)=\varphi(V^\top x) .
\]
Suppose that every function $F_V$ attains its minimum and that for some $\mathrm{gap}>0$
\begin{equation}\label{eq:hidden-gap}
 \ip{v_n}x=0\implies F_V(x)-\min F_V\ge\mathrm{gap}.
\end{equation}
Then for every deterministic method with $n-1$ oracle calls in $\R^d$ there exists $V$ such that $F_V(\hat x)-\min F_V\ge\mathrm{gap}$.
\item[(b)] Let a chain $\varphi\in\F_{L_0,L_1}^n$ have a minimizer of norm at most $R$, and let
\[
 \varphi(y)-\min\varphi\ge\mathrm{gap}>\varepsilon
\]
whenever $y_n=0$. Then $\mathcal N(\mathcal P_{\rm R})\ge n$. If, in addition, $\varphi(0)-\min\varphi\le\Delta$, then $\mathcal N(\mathcal P_\Delta)\ge n$.
\end{enumerate}
\end{proposition}

Part~(a) is the orthogonal invariance principle, that is, a reduction based on the orthogonal invariance of the class, in the form of the adversarial rotation argument of \citet[Section~3.3]{carmon2020lower}. The vectors $v_k$ are chosen during the run of the method orthogonally to the queries already made, so the oracle answers, including the function values, do not depend on the vectors not yet chosen. The quadratic chain for methods whose points lie in the linear span of the gradients is treated, for example, in \citet[Exercise~1.3]{gasnikov2021modern}. The transition to arbitrary deterministic methods goes back to \citet{nemirovski1983problem}. Part~(b) follows from~(a) with $d=2n$.

\subsection{Chain with geometrically decreasing slopes}\label{sec:ordered}

The links of this chain are affine functions of the coordinates with slopes $G_i$ decreasing geometrically. As a result, the shifts of the links are almost equal, and almost $L_1R$ links fit into a ball of radius $R$: only a factor $\sqrt{\log(L_1R)}$ is lost. What is new here is the choice of the function $H$, whose strong convexity modulus decreases exactly as the $(L_0,L_1)$-smoothness condition requires, and the choice of the slopes and shifts.

\paragraph{Slopes and shifts.} Fix $n\ge1$ and for $i=1,\dots,n$ set
\begin{equation}\label{eq:slopes}
\begin{aligned}
 G_i&=\varepsilon L_1\exp\Bigl((n-i)\sqrt{\frac{\log(4n)}n}\Bigr),\\
 \sigma_i&=\frac2{L_1}\log\Bigl(1+\frac{L_1G_1}{L_0}\Bigr)+\frac2{L_1}+\frac{4\log(1+\sqrt n)}{L_1G_i}\sum_{k=i}^nG_k .
\end{aligned}
\end{equation}
The slopes decrease from $G_1$ to $G_n=\varepsilon L_1$. By the definition of the shifts, for $j<l$
\begin{equation}\label{eq:shift-gap}
\begin{aligned}
 G_j\sigma_j-G_l\sigma_l
 &=(G_j-G_l)\Bigl[\frac2{L_1}\log\Bigl(1+\frac{L_1G_1}{L_0}\Bigr)+\frac2{L_1}\Bigr]+\frac{4\log(1+\sqrt n)}{L_1}\sum_{k=j}^{l-1}G_k\\
 &\ge(G_j-G_l)\,\frac2{L_1}\log\Bigl(1+\frac{L_1G_1}{L_0}\Bigr)+\frac{4\log(1+\sqrt n)}{L_1}\,G_j .
\end{aligned}
\end{equation}

\paragraph{Function on a polytope.} Define the polytope
\[
 K=\Bigl\{z\in\R^n:\ z\ge0,\ \sum_{i=1}^n\frac{z_i}{G_i}\le1\Bigr\}
\]
and, for $p\ge0$, the function
\[
 r(p)=\frac4{L_1}\Bigl[\Bigl(p+\frac1{\sqrt n}\Bigr)\log\bigl(1+\sqrt n\,p\bigr)-p\Bigr].
\]
Its derivatives are
\begin{equation}\label{eq:r-derivatives}
\begin{aligned}
 &r'(p)=\frac4{L_1}\log\bigl(1+\sqrt n\,p\bigr),\\
 &r''(p)=\frac{4\sqrt n}{L_1\bigl(1+\sqrt n\,p\bigr)},
\end{aligned}
\end{equation}
so that $r(0)=r'(0)=0$. For a point $z\in K$ denote $Z_n=z_1+\dots+z_n$ and set
\[
 H(z)=\begin{cases}
 2\Phi(Z_n)+Z_n\displaystyle\sum_{i=1}^nr\Bigl(\frac{z_i}{Z_n}\Bigr), & Z_n>0,\\[2pt]
 0, & z=0,
 \end{cases}
\]
where $\Phi$ is the function~\eqref{eq:Phi} from Section~\ref{sec:upper}. Finally, set
\begin{equation}\label{eq:ordered}
 \varphi(y)=\max_{z\in K}\bigl\{\ip{y+\sigma}z-H(z)\bigr\}.
\end{equation}

At the vertex $G_ie_i$ of the polytope, the expression under the maximum equals $G_i(y_i+\sigma_i)-H(G_ie_i)$, that is, the $i$-th link is an affine function with slope $G_i$.

The first term of $H$ depends only on the sum of the coordinates of $z$ and is responsible for the growth of the curvature together with the gradient norm. The second term is responsible for the distribution of this sum over the coordinates. The derivative of $r$ at zero is finite, so distant links are not activated at all.

\begin{lemma}[Membership in the class]\label{lem:ord-class}
The function $H$ satisfies the assumption of Lemma~\ref{lem:conj}. Consequently, $\varphi\in\F_{L_0,L_1}^n$ and its gradient satisfies inequality~\eqref{eq:conj-lip}.
\end{lemma}

\begin{proof}
\leavevmode
\begin{itemize}
\item \emph{Continuity.} The function $r$ is convex and $r(0)=0$, so on $[0,1]$
\[
 0\le r(p)\le p\,r'(p)\le p\,r'(1) .
\]
Hence, the second term of $H$ is at most $Z_nr'(1)$ and tends to zero together with $Z_n$. At the other points of $K$ continuity is obvious.
\item \emph{Second derivative along a segment.} Let $z,z'\in K$ be distinct, $u=z'-z$ and $w=z+tu$, $t\in[0,1]$. Denote $W_n=w_1+\dots+w_n$ and $U_n=u_1+\dots+u_n$, so that $W_n$ depends on $t$ and its derivative in $t$ equals $U_n$. Set
\[
 p=\frac w{W_n} .
\]
Coordinates that vanish at both points vanish on the whole segment, and the others are positive for $0<t<1$. Therefore, on $(0,1)$ the function $t\mapsto H(w)$ is twice continuously differentiable. The second term of $H$ is the perspective of the function $r(p_1)+\dots+r(p_n)$. Along the segment
\begin{align*}
 \frac{dp}{dt}&=\frac{u-pU_n}{W_n},\\
 \Phi''(W_n)&\by{=}{\eqref{eq:Phi-derivatives}}\frac1{L_0+L_1W_n},
\end{align*}
and differentiation gives
\begin{equation}\label{eq:H-second}
 \frac{d^2}{dt^2}H(w)=\frac{2U_n^2}{L_0+L_1W_n}+\sum_{i=1}^n\frac{r''(p_i)}{W_n}\,(u_i-p_iU_n)^2 .
\end{equation}
The vector $p$ lies in the simplex, so each of its coordinates is at most its norm and
\begin{equation}\label{eq:simplex-p}
 \frac1{\sqrt n}\le\norm p\le1 .
\end{equation}
Moreover, $\norm w=W_n\norm p$. We bound the coefficients in~\eqref{eq:H-second}:
\begin{align*}
 \frac{r''(p_i)}{W_n}
 &\by{=}{\eqref{eq:r-derivatives}}\frac{4\sqrt n}{L_1W_n\,(1+\sqrt n\,p_i)}
 \by{\ge}{\eqref{eq:simplex-p}}\frac2{L_1W_n\norm p}\\
 &\bypad{=}{}{\eqref{eq:r-derivatives}}\frac2{L_1\norm w}
 \ge\frac2{L_0+L_1\norm w},\\
 \frac2{L_0+L_1W_n}
 &\by{\ge}{\norm p\le1}\frac{2\norm p^2}{L_0+L_1W_n\norm p}
 =\frac{2\norm p^2}{L_0+L_1\norm w} .
\end{align*}
Substituting them into~\eqref{eq:H-second}, we obtain
\begin{align*}
 \frac{d^2}{dt^2}H(w)
 &\ge\frac2{L_0+L_1\norm w}\Bigl(\norm p^2U_n^2+\norm{u-pU_n}^2\Bigr)\\
 &\by{\ge}{(\mathrm a)}\frac{\norm u^2}{L_0+L_1\norm w}
 \by{\ge}{(\mathrm b)}\frac{\norm u^2}{L_0+L_1\max\{\norm z,\norm{z'}\}} .
\end{align*}
Here step (a) uses the inequality $\norm{u}^2\le2\norm{pU_n}^2+2\norm{u-pU_n}^2$, and step (b) uses the convexity of the norm.
\item \emph{Conclusion.} The restriction of $H$ to the segment minus the quadratic term from Lemma~\ref{lem:conj} has a nonnegative second derivative on $(0,1)$ and is continuous on $[0,1]$. Hence, it is convex, and Lemma~\ref{lem:conj} applies.\qedhere
\end{itemize}
\end{proof}

The second derivative of the term $2\Phi(Z_n)$ yields the factor $L_0+L_1\norm z$ in the strong convexity modulus and, consequently, the growth of the curvature of $\varphi$ together with the gradient norm. The factor $\sqrt n$ in the function $r$ matches the second term with the first one and determines the quantity
\[
 r'(1)=\frac4{L_1}\log\bigl(1+\sqrt n\bigr),
\]
on which the length of the chain depends.

\begin{lemma}[Chain]\label{lem:ord-local}
The function $\varphi$ is a chain: if $y\in\mathcal K_{j-1}$, then $\nabla\varphi(y)\in\mathcal K_j$.
\end{lemma}

\begin{proof}
Let $y\in\mathcal K_{j-1}$, that is, $y_i=0$ for $i\ge j$, and let $z=\nabla\varphi(y)$. We need to prove that $z_l=0$ for $l>j$. If $z=0$, there is nothing to prove. Otherwise $Z_n>0$. Let $p$ be the vector $z$ divided by $Z_n$.
\begin{itemize}
\item \emph{Optimality conditions.} In a neighborhood of $z$ the function $H$ is differentiable and $\partial_iH(z)=\nu+r'(p_i)$, where
\[
 \nu=2\Phi'(Z_n)+\sum_{k=1}^n\bigl[r(p_k)-p_kr'(p_k)\bigr].
\]
The point $z$ maximizes a concave function over the polytope $K$. By the Karush--Kuhn--Tucker theorem, there exists a multiplier $\lambda\ge0$ of the linear constraint defining $K$ such that
\begin{equation}\label{eq:ord-kkt}
 y_i+\sigma_i-\frac\lambda{G_i}\ \begin{cases}=\nu+r'(p_i), & z_i>0,\\ \le\nu, & z_i=0.\end{cases}
\end{equation}
\item \emph{Upper bound on $\nu$.} Since $r(p)\le p\,r'(p)$, we have $\nu\le2\Phi'(Z_n)$. Further, $Z_n\le G_1$, since $z\in K$ and all slopes are at most $G_1$. The function $\Phi'$ is increasing, therefore
\begin{equation}\label{eq:ord-nu-upper}
 \nu\le2\Phi'(G_1)\by{=}{\eqref{eq:Phi-derivatives}}\frac2{L_1}\log\Bigl(1+\frac{L_1G_1}{L_0}\Bigr).
\end{equation}
\item \emph{Lower bound on $\nu$.} Since $0\le r'(p_j)\le r'(1)$ and $y_j=0$, in both cases~\eqref{eq:ord-kkt} gives
\begin{equation}\label{eq:ord-nu-lower}
 -\lambda\le G_j\nu-G_j\sigma_j+G_jr'(1) .
\end{equation}
\item \emph{Coordinate $l$ is zero.} Since $y_l=0$ and $G_j\ge G_l$,
\begin{align*}
 G_l\Bigl(\sigma_l-\frac\lambda{G_l}-\nu\Bigr)
 &\bypad{\le}{\eqref{eq:ord-nu-lower}}{\eqref{eq:ord-nu-upper}}G_l\sigma_l-G_j\sigma_j+G_jr'(1)+(G_j-G_l)\nu\\
 &\by{\le}{\eqref{eq:ord-nu-upper}}G_l\sigma_l-G_j\sigma_j+G_jr'(1)+(G_j-G_l)\,\frac2{L_1}\log\Bigl(1+\frac{L_1G_1}{L_0}\Bigr)\\
 &\bypad{\le}{\eqref{eq:shift-gap}}{\eqref{eq:ord-nu-upper}}0 .
\end{align*}
If $z_l>0$, then by~\eqref{eq:ord-kkt} the left-hand side would equal $G_lr'(p_l)>0$. Hence, $z_l=0$.\qedhere
\end{itemize}
\end{proof}

\begin{lemma}\label{lem:ord-gap}
\leavevmode
\begin{enumerate}
\item[(a)] $\min\varphi=0$, and the minimum is attained at the point $-\sigma$.
\item[(b)] If $y_n=0$, then $\varphi(y)\ge2\varepsilon$.
\end{enumerate}
\end{lemma}

\begin{proof}
\leavevmode
\begin{itemize}
\item \emph{Statement~(a).} The point $z=0$ is feasible, so $\varphi\ge-H(0)=0$. On the other hand, $H\ge0$, and $\varphi(-\sigma)=\max_K\{-H\}=0$.
\item \emph{Statement~(b).} Substitute the vertex $G_ne_n$ into~\eqref{eq:ordered}. The function $\Phi$ is convex, $\Phi(0)=0$ and $\Phi'$ is increasing, so $\Phi(G_n)\le G_n\Phi'(G_n)\le G_n\Phi'(G_1)$. Moreover, $r(1)\le r'(1)$. Hence,
\begin{align*}
 \varphi(y)
 &\ge G_n\sigma_n-2\Phi(G_n)-G_nr(1)\\
 &\ge G_n\sigma_n-\frac{2G_n}{L_1}\log\Bigl(1+\frac{L_1G_1}{L_0}\Bigr)-G_n\,\frac4{L_1}\log\bigl(1+\sqrt n\bigr)
 \by{=}{\eqref{eq:slopes}}\frac{2G_n}{L_1}=2\varepsilon .\qedhere
\end{align*}
\end{itemize}
\end{proof}

\begin{lemma}[Distance to the minimizer]\label{lem:ord-radius}
\[
 L_1\norm\sigma\le2\sqrt n\,\log\Bigl(1+\frac{L_1^2\varepsilon}{L_0}\Bigr)+9\,n\sqrt{\log(4n)} .
\]
\end{lemma}

\begin{proof}
We bound the three quantities that make up the shifts $\sigma_i$:
\begin{align*}
 &2\log\Bigl(1+\frac{L_1G_1}{L_0}\Bigr)\by{\le}{(\mathrm a)}2\log\Bigl(1+\frac{L_1^2\varepsilon}{L_0}\Bigr)+2\sqrt{n\log(4n)},\\
 &4\log\bigl(1+\sqrt n\bigr)\by{\le}{(\mathrm b)}2\log(4n),\\
 &\sum_{k=i}^n\frac{G_k}{G_i}\by{\le}{(\mathrm c)}1+\sqrt{\frac n{\log(4n)}} .
\end{align*}
Here step (a) uses the inequality $\log(1+xe^v)\le\log(1+x)+v$ for $x,v\ge0$, step (b) uses the inequality $(1+\sqrt n)^2\le4n$, and step (c) uses the sum of a geometric progression together with the inequality $e^{-x}(1+x)\le1$. Adding these bounds, we obtain for every coordinate
\[
 L_1\sigma_i\le2\log\Bigl(1+\frac{L_1^2\varepsilon}{L_0}\Bigr)+4\sqrt{n\log(4n)}+2\log(4n)+2 .
\]
All $\sigma_i$ are positive, so $\norm\sigma\le\sqrt n\max_i\sigma_i$ and
\[
 L_1\norm\sigma\le2\sqrt n\log\Bigl(1+\frac{L_1^2\varepsilon}{L_0}\Bigr)+n\sqrt{\log(4n)}\,\Bigl[4+2\sqrt{\frac{\log(4n)}n}+\frac2{\sqrt{n\log(4n)}}\Bigr].
\]
The ratio of $\log(4n)$ to $n$ decreases and is therefore at most $\log4<2$, while $n\log(4n)\ge\log4>1$. Hence, the expression in brackets is less than $4+2\sqrt2+2<9$.
\end{proof}

\begin{theorem}[Lower bound for the chain with decreasing slopes]\label{thm:ordered}
Let a positive integer $n$ satisfy
\begin{equation}\label{eq:ordered-condition}
 2\sqrt n\,\log\Bigl(1+\frac{L_1^2\varepsilon}{L_0}\Bigr)+9\,n\sqrt{\log(4n)}\le L_1R .
\end{equation}
Then for every deterministic method with $n-1$ oracle calls in $\R^{2n}$ there exists a function $f\in\mathcal P_{\rm R}$ on $\R^{2n}$ with $\min f=0$ and $f(\hat x)\ge2\varepsilon$. In particular, $\mathcal N(\mathcal P_{\rm R})\ge n$.
\end{theorem}

\begin{proof}
The function $\varphi$ is a convex chain from the class $\F_{L_0,L_1}^n$ (Lemmas~\ref{lem:ord-class} and~\ref{lem:ord-local}). Its minimum is zero and is attained at the point $-\sigma$, and for $y_n=0$ its value is at least $2\varepsilon$ (Lemma~\ref{lem:ord-gap}). Lemma~\ref{lem:ord-radius} and condition~\eqref{eq:ordered-condition} imply that the norm of $\sigma$ is at most $R$. It remains to apply Proposition~\ref{prop:reduction}(b) with $\mathrm{gap}=2\varepsilon$.
\end{proof}

The initial optimality gap of these functions is not bounded by $4\varepsilon$. Substituting the vertex $G_1e_1$ into~\eqref{eq:ordered}, as in the proof of Lemma~\ref{lem:ord-gap}(b), gives
\[
 \varphi(0)\ge\frac{2G_1}{L_1}=2\varepsilon\exp\Bigl((n-1)\sqrt{\frac{\log(4n)}n}\Bigr),
\]
that is, the gap is exponentially large in $\sqrt{n\log n}$. Therefore, Theorem~\ref{thm:ordered} does not contradict the upper bound with the term $(L_1R)^{2/3}$ on the class $\mathcal P_{4\varepsilon}$ (Appendix~\ref{app:gap}).

\subsection{Main lower bound}\label{sec:lower-radius-proof}

It remains to choose a chain length satisfying condition~\eqref{eq:ordered-condition}.

\begin{lemma}[Choice of the chain length]\label{lem:choice}
If $L_0R^2\ge\varepsilon$ and $L_1R\ge16$, then
\[
 1+\mathcal N(\mathcal P_{\rm R})\ge\frac{L_1R}{24\sqrt{\log(2+L_1R)}} .
\]
\end{lemma}

\begin{proof}
Set
\[
 n=\Bigl\lfloor\frac{L_1R}{24\sqrt{\log(2+L_1R)}}\Bigr\rfloor .
\]
If $n=0$, the statement is obvious. Let $n\ge1$. We verify condition~\eqref{eq:ordered-condition} by bounding its two terms.
\begin{itemize}
\item \emph{Second term.} The number $4n$ is less than $2+L_1R$, so $\log(4n)\le\log(2+L_1R)$ and
\[
 9n\sqrt{\log(4n)}\le9n\sqrt{\log(2+L_1R)}\by{\le}{\text{def. of }n}\frac38L_1R .
\]
\item \emph{First term.} The ratio of $\log^{3/2}(2+x)$ to $x$ decreases for $x\ge16$ and is less than three eighths at the point $16$. Therefore,
\begin{align*}
 2\sqrt n\log\Bigl(1+\frac{L_1^2\varepsilon}{L_0}\Bigr)
 &\by{\le}{L_0R^2\ge\varepsilon}2\sqrt n\log\bigl(1+(L_1R)^2\bigr)
 \by{\le}{1+t^2\le(2+t)^2}4\sqrt n\log(2+L_1R)\\
 &\by{\le}{\text{def. of }n}\frac4{\sqrt{24}}\,L_1R\,\sqrt{\frac{\log^{3/2}(2+L_1R)}{L_1R}}\\
 &\by{<}{L_1R\ge16}\frac4{\sqrt{24}}\sqrt{\frac38}\,L_1R=\frac{L_1R}2 .
\end{align*}
\end{itemize}
The sum of the terms is less than $L_1R$, and Theorem~\ref{thm:ordered} gives
\[
 1+\mathcal N(\mathcal P_{\rm R})\ge1+n\ge\frac{L_1R}{24\sqrt{\log(2+L_1R)}} .\qedhere
\]
\end{proof}

\begin{theorem}[Lower bound for the radius-only class]\label{thm:lower-radius}
If $L_0R^2\ge\varepsilon$, then
\begin{equation}\label{eq:lower-radius}
 1+\mathcal N(\mathcal P_{\rm R})\ge\frac1{29}\Bigl[1+\sqrt{\frac{L_0R^2}\varepsilon}+\frac{L_1R}{\sqrt{\log(2+L_1R)}}\Bigr].
\end{equation}
\end{theorem}

\begin{proof}
We bound each of the three terms in brackets by $1+\mathcal N(\mathcal P_{\rm R})$.
\begin{itemize}
\item \emph{First.} Obviously, $1\le1+\mathcal N(\mathcal P_{\rm R})$.
\item \emph{Second.} Functions with $L_0$-Lipschitz gradient belong to the class. The classical bound \citep[Theorem~2.1.7]{nesterov2018lectures} together with Proposition~\ref{prop:reduction} gives
\[
 \sqrt{\frac{L_0R^2}\varepsilon}\le4\bigl(1+\mathcal N(\mathcal P_{\rm R})\bigr).
\]
\item \emph{Third.} For $L_1R\ge16$, by Lemma~\ref{lem:choice} it is at most $24(1+\mathcal N(\mathcal P_{\rm R}))$. For $L_1R<16$ it is less than $10$, since it increases in $L_1R$ and
\[
 \frac{16}{\sqrt{\log18}}\by{<}{\log18>4\log2}\frac{16}{\sqrt{4\log2}}\by{<}{\log2>2/3}\sqrt{96}<10 .
\]
\end{itemize}
Adding up, we obtain that the expression in brackets in~\eqref{eq:lower-radius} is at most $29(1+\mathcal N(\mathcal P_{\rm R}))$.
\end{proof}

Together with Theorem~\ref{thm:upper}, for $L_0R^2\ge\varepsilon$ we obtain
\[
 1+\mathcal N(\mathcal P_{\rm R})=\widetilde\Theta\Bigl(1+\sqrt{\frac{L_0R^2}\varepsilon}+L_1R\Bigr).
\]
The exponents in the lower and upper bounds coincide, and the bounds differ only by logarithmic factors in the $L_1R$ term.

\section{Conclusion and open questions}\label{sec:conclusion}

Large gradients and the curvature growing together with them create a regime in which the number of oracle calls in the worst case depends on $L_1R$ and does not depend on $\varepsilon$. If, besides $L_0$, $L_1$ and $\varepsilon$, the method knows only an upper bound $R$ on the initial distance to a minimizer, then for $L_0R^2\ge\varepsilon$ the optimal complexity $N$ among deterministic first-order methods satisfies
\[
 1+N=\widetilde\Theta\Bigl(1+\sqrt{\frac{L_0R^2}\varepsilon}+L_1R\Bigr).
\]
In other words, the cost of passing through the region of large gradients is $L_1R$ up to logarithmic factors. If the initial optimality gap is at most $4\varepsilon$, this cost drops to $(L_1R)^{2/3}$ (Appendix~\ref{app:gap}). The following questions remain open.

\begin{itemize}
\item \emph{Logarithmic factors.} Between the lower and upper bounds there remains a factor of order $\log^{3/2}(2+L_1R)$ in the $L_1R$ term. The factor $\sqrt{\log(2+L_1R)}$ is lost in the lower bound because of the choice of the slopes and shifts of the chain, and the factor $\log(2+L_1R)$ in the upper bound because of the line search at the first stage. We do not know whether the large-gradient regime can be passed in $O(1+L_1R)$ calls without line search.
\item \emph{Known bound on the optimality gap.} For the class $\mathcal P_\Delta$ with $\Delta$ greater than $4\varepsilon$ but smaller than the right-hand side of~\eqref{eq:radial-gap}, only the upper bound of Theorem~\ref{thm:upper-gap}(b) is proved. A matching lower bound in terms of the dependence on $\Delta$ is unknown. The exact dependence of the complexity on $\Delta/\varepsilon$ for $\varepsilon<\Delta<4\varepsilon$ is not established either.
\item \emph{Dimension.} The lower bounds for arbitrary methods are realized in a dimension proportional to the number of calls, while in fixed dimension with fixed $L_0$, $L_1$ and $\varepsilon$ the complexity grows only logarithmically in $R$ (Appendix~\ref{app:fixed-dim}). The intermediate regime has not been studied.
\item \emph{Randomized methods and inexact gradients.} All lower bounds are obtained for deterministic methods and an exact oracle.
\end{itemize}

\clearpage
\bibliographystyle{plainnat}
\bibliography{refs}

\clearpage
\appendix
\section{Auxiliary inequalities}\label{app:tools}

In this appendix we prove Lemma~\ref{lem:class} and the auxiliary lemmas of Section~\ref{sec:upper}. In Subsection~\ref{app:tools-class} the function $\ell$ is arbitrary, as in Definition~\ref{def:class}. Starting from Subsection~\ref{app:tools-gronwall}, $f\in\F_{L_0,L_1}^d$ and $x^\ast$ is a minimizer of $f$.

\subsection{Characterization of the class}\label{app:tools-class}

\begin{proof}[Proof of Lemma~\ref{lem:class}\proofup{lem:class}]
\leavevmode
\begin{itemize}
\item \emph{Necessity in~(a).} The gradient of a convex differentiable function is continuous \citep[Theorem~25.5]{rockafellar1970convex}, so on every ball it is bounded by some number $M$. If the segment $[x,y]$ lies in such a ball, then by~\eqref{eq:class}
\[
 \norm{\nabla f(x)-\nabla f(y)}\le\ell(M)\norm{x-y} .
\]
Hence, the map $t\mapsto\nabla f(x+tv)$ is Lipschitz on every segment, and therefore locally absolutely continuous and differentiable for almost all $t$. At such a point $t$, apply~\eqref{eq:class} to the points $x+tv$ and $x+(t+h)v$, divide by $|h|$ and let $h$ tend to zero. By the continuity of $\nabla f$ and $\ell$ we obtain~\eqref{eq:line-derivative}.
\item \emph{Sufficiency in~(a).} Let $v=y-x$. By the fundamental theorem of calculus and the monotonicity of $\ell$
\begin{align*}
 \norm{\nabla f(y)-\nabla f(x)}
 &\le\int_0^1\ell\bigl(\norm{\nabla f(x+tv)}\bigr)\norm v\,dt\\
 &\le\ell\Bigl(\sup_{u\in[x,y]}\norm{\nabla f(u)}\Bigr)\norm{x-y} .
\end{align*}
\item \emph{Statement~(b).} For $f\in C^2$ the derivative of the map $t\mapsto\nabla f(x+tv)$ equals $\nabla^2f(x+tv)v$ and is continuous in $t$. Therefore, condition~\eqref{eq:line-derivative} for almost all $t$ is equivalent to the inequality
\[
 \norm{\nabla^2f(x)v}\le\ell\bigl(\norm{\nabla f(x)}\bigr)\norm v
\]
for all $x$ and $v$. The Hessian of a convex $C^2$ function is a symmetric positive semidefinite matrix, and for it this inequality is equivalent to~\eqref{eq:hessian-bound}. It remains to apply~(a).\qedhere
\end{itemize}
\end{proof}

\subsection{Rate of change of the gradient}\label{app:tools-gronwall}

The following lemma shows that the local smoothness constant $L_0+L_1\norm{\nabla f(y)}$ changes by a factor of at most $e^{L_1\norm{z-y}}$ when passing from the point $y$ to the point $z$. The remaining bounds of the appendix are derived from this fact.

\begin{lemma}[Consequences of Grönwall's inequality]\label{lem:gronwall}
For all points $y,z\in\R^d$
\begin{align*}
 &(\text{a})&L_0+L_1\norm{\nabla f(z)}&\le\bigl(L_0+L_1\norm{\nabla f(y)}\bigr)e^{L_1\norm{z-y}},\\
 &(\text{b})&\norm{\nabla f(z)-\nabla f(y)}&\le\bigl(L_0+L_1\norm{\nabla f(y)}\bigr)\frac{e^{L_1\norm{z-y}}-1}{L_1},\\
 &(\text{c})&0\le f(z)-f(y)-\ip{\nabla f(y)}{z-y}&\le\bigl(L_0+L_1\norm{\nabla f(y)}\bigr)\frac{e^{L_1\norm{z-y}}-1-L_1\norm{z-y}}{L_1^2} .
\end{align*}
\end{lemma}

\begin{proof}
We track the gradient at the points $y+t(z-y)$, $t\in[0,1]$.
\begin{itemize}
\item \emph{Integral inequality.} By Lemma~\ref{lem:class}(a), the map $t\mapsto\nabla f(y+t(z-y))$ is absolutely continuous, and for almost all $t$ the norm of its derivative is at most $(L_0+L_1\norm{\nabla f(y+t(z-y))})\norm{z-y}$. Therefore, for all $t\in[0,1]$
\begin{equation}\label{eq:gronwall-int}
 \norm{\nabla f\bigl(y+t(z-y)\bigr)-\nabla f(y)}\le\norm{z-y}\int_0^t\Bigl(L_0+L_1\norm{\nabla f\bigl(y+\tau(z-y)\bigr)}\Bigr)d\tau .
\end{equation}
By the triangle inequality this implies
\begin{align*}
 \MoveEqLeft L_0+L_1\norm{\nabla f\bigl(y+t(z-y)\bigr)}\\
 &\le L_0+L_1\norm{\nabla f(y)}+L_1\norm{z-y}\int_0^t\Bigl(L_0+L_1\norm{\nabla f\bigl(y+\tau(z-y)\bigr)}\Bigr)d\tau .
\end{align*}
\item \emph{Grönwall's inequality.} The left-hand side of the last inequality is continuous in $t$, and by Grönwall's inequality in integral form
\[
 L_0+L_1\norm{\nabla f\bigl(y+t(z-y)\bigr)}\le\bigl(L_0+L_1\norm{\nabla f(y)}\bigr)e^{L_1\norm{z-y}t} .
\]
For $t=1$ this is statement~(a).
\item \emph{Increment of the gradient.} Substitute this bound into~\eqref{eq:gronwall-int}:
\begin{align*}
 \norm{\nabla f\bigl(y+t(z-y)\bigr)-\nabla f(y)}
 &\le\bigl(L_0+L_1\norm{\nabla f(y)}\bigr)\norm{z-y}\int_0^te^{L_1\norm{z-y}\tau}d\tau\\
 &=\bigl(L_0+L_1\norm{\nabla f(y)}\bigr)\frac{e^{L_1\norm{z-y}t}-1}{L_1} .
\end{align*}
For $t=1$ this is statement~(b).
\item \emph{Increment of the function.} The left inequality in~(c) is the convexity of $f$. For the right one
\begin{align*}
 f(z)-f(y)-\ip{\nabla f(y)}{z-y}
 &=\int_0^1\ip{\nabla f\bigl(y+t(z-y)\bigr)-\nabla f(y)}{z-y}\,dt\\
 &\by{\le}{(\text{b})}\bigl(L_0+L_1\norm{\nabla f(y)}\bigr)\norm{z-y}\int_0^1\frac{e^{L_1\norm{z-y}t}-1}{L_1}\,dt\\
 &=\bigl(L_0+L_1\norm{\nabla f(y)}\bigr)\frac{e^{L_1\norm{z-y}}-1-L_1\norm{z-y}}{L_1^2} .
\end{align*}
The second step uses the Cauchy--Schwarz inequality and statement~(b).\qedhere
\end{itemize}
\end{proof}

For twice continuously differentiable functions, inequalities~(b) and~(c) were obtained by \citet[Lemma~2.5]{vankov2025optimizing}. The inequality $e^x-1-x\le x^2e^x/2$ for $x\ge0$ and statement~(c) imply
\begin{equation}\label{eq:local-descent}
 f(z)\le f(y)+\ip{\nabla f(y)}{z-y}+\frac{\bigl(L_0+L_1\norm{\nabla f(y)}\bigr)e^{L_1\norm{z-y}}}2\,\norm{z-y}^2 .
\end{equation}
At distances of order $1/L_1$, this is the usual inequality for functions with Lipschitz gradient \citep[Lemma~1.2.3]{nesterov2018lectures} with a constant proportional to $L_0+L_1\norm{\nabla f(y)}$.

\subsection{Descent lemma}\label{app:tools-descent}

\begin{proof}[Proof of Lemma~\ref{lem:safe-step}\proofup{lem:safe-step}]
\leavevmode
\begin{itemize}
\item \emph{Decrease of the function.} Since $L_0+L_1\norm{\nabla f(y)}\ge L_1\norm{\nabla f(y)}$, the step length
\[
 \norm{y^+-y}=\frac{\norm{\nabla f(y)}}{4\bigl(L_0+L_1\norm{\nabla f(y)}\bigr)}
\]
is at most $1/(4L_1)$, and the factor $e^{L_1\norm{y^+-y}}$ in~\eqref{eq:local-descent} is at most $e^{1/4}$. Therefore,
\begin{align*}
 f(y^+)
 &\by{\le}{\eqref{eq:local-descent}}f(y)-\frac{\norm{\nabla f(y)}^2}{4\bigl(L_0+L_1\norm{\nabla f(y)}\bigr)}+\frac{e^{1/4}\norm{\nabla f(y)}^2}{32\bigl(L_0+L_1\norm{\nabla f(y)}\bigr)}\\
 &\bypad{=}{}{\eqref{eq:local-descent}}f(y)-\Bigl(\frac14-\frac{e^{1/4}}{32}\Bigr)\frac{\norm{\nabla f(y)}^2}{L_0+L_1\norm{\nabla f(y)}}\\
 &\bypad{\le}{(\mathrm a)}{\eqref{eq:local-descent}}f(y)-\frac{5\norm{\nabla f(y)}^2}{24\bigl(L_0+L_1\norm{\nabla f(y)}\bigr)} .
\end{align*}
Here step (a) uses the inequality $e^{1/4}<4/3$, which follows from $e^t\le1/(1-t)$ for $0\le t<1$.
\item \emph{Decrease of the distance to the minimizer.} We have
\begin{align*}
 \ip{\nabla f(y)}{y-x^\ast}
 &\by{\ge}{(\mathrm b)}f(y)-f^\ast
 \by{\ge}{(\mathrm c)}f(y)-f(y^+)\\
 &\by{\ge}{(\mathrm d)}\frac{5\norm{\nabla f(y)}^2}{24\bigl(L_0+L_1\norm{\nabla f(y)}\bigr)} ,
\end{align*}
where step (b) uses convexity, step (c) uses the inequality $f(y^+)\ge f^\ast$, and step (d) uses the decrease of the function proved above.
Expanding the squared norm and substituting this bound, we obtain
\begin{align*}
 \norm{y^+-x^\ast}^2
 &=\norm{y-x^\ast}^2-\frac{\ip{\nabla f(y)}{y-x^\ast}}{2\bigl(L_0+L_1\norm{\nabla f(y)}\bigr)}+\frac{\norm{\nabla f(y)}^2}{16\bigl(L_0+L_1\norm{\nabla f(y)}\bigr)^2}\\
 &\le\norm{y-x^\ast}^2-\frac{5\norm{\nabla f(y)}^2}{48\bigl(L_0+L_1\norm{\nabla f(y)}\bigr)^2}+\frac{3\norm{\nabla f(y)}^2}{48\bigl(L_0+L_1\norm{\nabla f(y)}\bigr)^2} .\qedhere
\end{align*}
\end{itemize}
\end{proof}

\subsection{Bounds in a neighborhood of the minimizer}\label{app:tools-anchored}

The gradient vanishes at the minimizer, so Lemma~\ref{lem:gronwall} with $y=x^\ast$ bounds the growth of the function away from the minimizer.

\begin{lemma}[Growth around the minimizer]\label{lem:radial}
For every point $x\in\R^d$
\begin{align*}
 \norm{\nabla f(x)}&\le\frac{L_0}{L_1}\bigl(e^{L_1\norm{x-x^\ast}}-1\bigr),\\
 f(x)-f^\ast&\le\frac{L_0}{L_1^2}\bigl(e^{L_1\norm{x-x^\ast}}-1-L_1\norm{x-x^\ast}\bigr).
\end{align*}
In particular, if $\norm{x^0-x^\ast}\le R$, then
\begin{equation}\label{eq:radial-gap}
 f(x^0)-f^\ast\le\frac{L_0}{L_1^2}\bigl(e^{L_1R}-1-L_1R\bigr).
\end{equation}
\end{lemma}

\begin{proof}
Both bounds are statements~(b) and~(c) of Lemma~\ref{lem:gronwall} with $y=x^\ast$ and $z=x$, since $\nabla f(x^\ast)=0$. The right-hand side of the second bound increases in $\norm{x-x^\ast}$, which gives~\eqref{eq:radial-gap}.
\end{proof}

Bound~\eqref{eq:radial-gap} is tight. The function
\[
 f(x)=\frac{L_0}{L_1^2}\bigl(e^{L_1\norm x}-1-L_1\norm x\bigr)
\]
is twice continuously differentiable, and
\[
 \norm{\nabla f(x)}=\frac{L_0}{L_1}\bigl(e^{L_1\norm x}-1\bigr).
\]
For $x\ne0$ its curvature along the radius equals
\[
 L_0e^{L_1\norm x}=L_0+L_1\norm{\nabla f(x)} ,
\]
and in the tangential directions the curvature equals
\[
 \frac{\norm{\nabla f(x)}}{\norm x}=\frac{L_0\bigl(e^{L_1\norm x}-1\bigr)}{L_1\norm x}\le L_0e^{L_1\norm x} ,
\]
since $e^s-1\le se^s$. At zero the Hessian equals $L_0I$, and both curvatures tend to $L_0$. By Lemma~\ref{lem:class}(b) the function belongs to the class $\F_{L_0,L_1}^d$, and for $\norm{x^0}=R$ equality holds in~\eqref{eq:radial-gap}. Thus, the initial optimality gap may be exponentially large in $L_1R$.

\begin{proof}[Proof of Lemma~\ref{lem:anchored}\proofup{lem:anchored}]
For $\nabla f(y)=0$ the statement coincides with the convexity inequality. Let $\nabla f(y)\ne0$. We compare the affine approximation at $y$ with the growth bound of Lemma~\ref{lem:radial} at the points of the ray
\[
 z_t=x^\ast+\frac{t\,\nabla f(y)}{\norm{\nabla f(y)}},
\]
where $t\ge0$.
\begin{itemize}
\item \emph{Bound for fixed $t$.} Since $\ip{\nabla f(y)}{z_t-x^\ast}=t\norm{\nabla f(y)}$,
\begin{align*}
 f(y)+\ip{\nabla f(y)}{x^\ast-y}+t\norm{\nabla f(y)}
 &\bypad{=}{}{(\mathrm a)}f(y)+\ip{\nabla f(y)}{z_t-y}\\
 &\by{\le}{(\mathrm a)}f(z_t)
 \by{\le}{(\mathrm b)}f^\ast+\frac{L_0}{L_1^2}\bigl(e^{L_1t}-1-L_1t\bigr),
\end{align*}
where step (a) uses convexity and step (b) uses Lemma~\ref{lem:radial}.
\item \emph{Choice of $t$.} Set
\[
 t=\frac1{L_1}\log\Bigl(1+\frac{L_1\norm{\nabla f(y)}}{L_0}\Bigr).
\]
Then
\[
 e^{L_1t}=1+\frac{L_1\norm{\nabla f(y)}}{L_0},
\]
and a direct computation using formula~\eqref{eq:Phi} gives
\[
 t\norm{\nabla f(y)}-\frac{L_0}{L_1^2}\bigl(e^{L_1t}-1-L_1t\bigr)=\Phi\bigl(\norm{\nabla f(y)}\bigr).
\]
Substituting this into the previous bound, we obtain the statement of the lemma.\qedhere
\end{itemize}
\end{proof}

The chosen $t$ maximizes the left-hand side of the last equality, that is, $\Phi$ is the conjugate of the growth function from Lemma~\ref{lem:radial}. The inequality of Lemma~\ref{lem:anchored} is claimed only at the point $x^\ast$: the affine function with the additional term $\Phi$ need not be a global minorant of $f$.

\subsection{Inexact line search}\label{app:tools-search}

\begin{proof}[Proof of Lemma~\ref{lem:search}\proofup{lem:search}]
Denote by $z(t)=v+t(x-v)$ the points of the segment $[v,x]$. The function $t\mapsto f(z(t))$ is convex on $[0,1]$, and its derivative $\ip{\nabla f(z(t))}{x-v}$ is nondecreasing. For the point $y=z(t)$
\[
 \ip{\nabla f(y)}{y-v}=t\,\ip{\nabla f(z(t))}{x-v} .
\]
\begin{itemize}
\item \emph{Stopping criterion.} If it was triggered at a point $z\in[v,x]$, then $z$ lies in the ball $B(x^\ast,r)$, and by convexity
\[
 f(z)-f^\ast\le\ip{\nabla f(z)}{z-x^\ast}\le r\norm{\nabla f(z)}\le\delta .
\]
From now on we assume that the stopping criterion was not triggered.
\item \emph{Return at the endpoints of the segment.} If $\ip{\nabla f(v)}{x-v}\ge0$, then the derivative is nonnegative on $[0,1]$. Hence, $f(v)\le f(x)$, and for $y=v$ the left-hand side of the second inequality of~\eqref{eq:search} is zero. If $\ip{\nabla f(x)}{x-v}\le0$, then for $y=x$ the first inequality of~\eqref{eq:search} is obvious, and the left-hand side of the second one is nonpositive.
\item \emph{Bisection.} The algorithm maintains the inequalities
\[
 \ip{\nabla f(z(t_-))}{x-v}<0<\ip{\nabla f(z(t_+))}{x-v} .
\]
If the derivative at the midpoint is zero, then this midpoint minimizes the function on the segment, and both inequalities~\eqref{eq:search} hold for it. Otherwise, after $m$ steps $t_+-t_-=2^{-m}$.
\item \emph{Decrease for $y=z(t_+)$.} The derivative at the point $t_+$ is positive, so the function is nondecreasing on $[t_+,1]$, and $f(y)\le f(x)$.
\item \emph{Error for $y=z(t_+)$.} The distance between the points $z(t_+)$ and $z(t_-)$ equals $2^{-m}\norm{x-v}$, and $\norm{x-v}\le2r$. By condition~\eqref{eq:search-m}
\[
 L_1\norm{z(t_+)-z(t_-)}\le2^{1-m}L_1r\le\frac{2L_1r}{128\bigl(1+L_1r\bigr)^2}\le\frac1{32} .
\]
By Lemma~\ref{lem:gronwall}(b) with base point $z(t_+)$ and the inequality $e^s-1\le se^s$
\[
 \norm{\nabla f(z(t_+))-\nabla f(z(t_-))}\le2\bigl(L_0+L_1\norm{\nabla f(y)}\bigr)\norm{z(t_+)-z(t_-)} .
\]
Since $\ip{\nabla f(z(t_-))}{x-v}<0$ and $t_+\le1$,
\begin{align*}
 \ip{\nabla f(y)}{y-v}
 &\bypad{\le}{}{\eqref{eq:search-m}}t_+\ip{\nabla f(z(t_+))-\nabla f(z(t_-))}{x-v}\\
 &\bypad{\le}{(\mathrm a)}{\eqref{eq:search-m}}2r\norm{\nabla f(z(t_+))-\nabla f(z(t_-))}\\
 &\bypad{\le}{}{\eqref{eq:search-m}}2^{3-m}r^2\bigl(L_0+L_1\norm{\nabla f(y)}\bigr)\\
 &\by{\le}{\eqref{eq:search-m}}\frac{r^2\bigl(L_0+L_1\norm{\nabla f(y)}\bigr)}{16\bigl(1+L_0r^2/\delta+L_1r\bigr)^2} ,
\end{align*}
where step (a) uses the inequalities $t_+\le1$ and $\norm{x-v}\le2r$.
\item \emph{Relation to the gradient.} The stopping criterion does not hold at $y$, that is, $r\norm{\nabla f(y)}>\delta$. Therefore,
\begin{align*}
 \frac{r\bigl(L_0+L_1\norm{\nabla f(y)}\bigr)}{\norm{\nabla f(y)}}
 &\bypad{=}{}{(\mathrm b)}\frac{L_0r}{\norm{\nabla f(y)}}+L_1r\\
 &\by{<}{(\mathrm b)}\frac{L_0r^2}\delta+L_1r
 <1+\frac{L_0r^2}\delta+L_1r ,
\end{align*}
where step (b) uses the inequality $r\norm{\nabla f(y)}>\delta$. Squaring and dividing by $16(L_0+L_1\norm{\nabla f(y)})$, we obtain
\[
 \frac{r^2\bigl(L_0+L_1\norm{\nabla f(y)}\bigr)}{16\bigl(1+L_0r^2/\delta+L_1r\bigr)^2}<\frac{\norm{\nabla f(y)}^2}{16\bigl(L_0+L_1\norm{\nabla f(y)}\bigr)} ,
\]
and the second inequality of~\eqref{eq:search} holds.
\item \emph{Number of calls.} The algorithm makes two calls at the endpoints of the segment and at most $m$ calls at midpoints.\qedhere
\end{itemize}
\end{proof}

\section{Class with bounded initial optimality gap}\label{app:gap}

In the main text, only an upper bound $R$ on the distance from the starting point to a minimizer is known, and the initial optimality gap may be exponentially large in $L_1R$. In this appendix we consider the class $\mathcal P_\Delta$ from Section~\ref{sec:setting}, in which the initial optimality gap is bounded by a known number $\Delta$, and primarily the case $\Delta=4\varepsilon$. The exponent of $L_1R$ in the complexity bound then drops from one to $2/3$.

The condition $f(x^0)-f^\ast\le4\varepsilon$ is natural in restart schemes: each phase starts from a point where the optimality gap is already bounded by a constant multiple of the accuracy required from this phase.

The lower bound is proved in Subsection~\ref{sec:eq-chain} on a chain with equal slopes of the links. The upper bounds are obtained by a variant of the accelerated method of Section~\ref{sec:upper} with other ways of counting iterations (Subsections~\ref{sec:gap-method}--\ref{sec:count2}). The final statements are formulated at the end of the corresponding subsections.

\subsection{Lower bound: chain with equal slopes}\label{sec:eq-chain}

As in Subsection~\ref{sec:ordered}, we construct a chain by conjugation on a polytope and apply Lemma~\ref{lem:conj}. Here the slopes of all links are equal, and the shifts form an arithmetic progression with common difference $-b$.

\paragraph{Slope and shifts.} Fix $n\ge1$ and set
\begin{equation}\label{eq:eq-chain-params}
\begin{aligned}
 &b=R\Bigl(\sum_{k=n+1}^{2n}k^2\Bigr)^{-1/2},\\
 &\sigma_i=(2n-i)\,b,\\
 &G=\frac{2\varepsilon}{nb} .
\end{aligned}
\end{equation}
where $i=1,\dots,n$. The sum in the definition of $b$ lies between $n^3$ and $4n^3$.

\paragraph{Function on a polytope.} Let a function $r\in C^1([0,1])$ be strongly convex with some positive modulus, $r(0)=r'(0)=0$ and $r'(1)=b$. Consider the polytope
\[
 K=\bigl\{z\in\R^n:\ z\ge0,\ Z_n\le G\bigr\},
\]
where $Z_n=z_1+\dots+z_n$. For $z\in K$ set
\[
\begin{aligned}
 p_i&=\frac{z_i}G,\\
 p_0&=1-\frac{Z_n}G ,
\end{aligned}
\]
where $i=1,\dots,n$.
The coordinate $p_0$ is the normalized slack in the constraint $Z_n\le G$, and it allows $\psi$ to have a whole region of minimizers. Define a function on $K$ and its conjugate:
\begin{equation}\label{eq:eq-chain}
\begin{aligned}
 &H(z)=G\Bigl[\sum_{i=0}^nr(p_i)-r(1)\Bigr],\\
 &\psi(y)=\max_{z\in K}\bigl\{\ip{y+\sigma}z-H(z)\bigr\}.
\end{aligned}
\end{equation}
Here $H(0)=H(Ge_i)=0$.

\begin{lemma}[Properties of the chain with equal slopes]\label{lem:eq-chain}
The function $\psi$ is convex, continuously differentiable, and
\begin{enumerate}
\item[(a)] its minimum is zero, and the set of minimizers is
\[
 \bigl\{y\in\R^n:\ y_i\le-(\sigma_i+b)\ \text{for}\ i=1,\dots,n\bigr\}.
\]
The minimizer closest to zero, $-b\,(2n,2n-1,\dots,n+1)$, has norm $R$,
\item[(b)] $\psi(0)=(4-2/n)\varepsilon$,
\item[(c)] $\psi(y)\ge2\varepsilon$ whenever $y_n=0$,
\item[(d)] $\psi$ is a chain,
\item[(e)] $\nabla\psi(y)=z(y)=G\,(p_1,\dots,p_n)$, where $z(y)$ is the unique maximizer in~\eqref{eq:eq-chain}.
\end{enumerate}
\end{lemma}

\begin{proof}
\leavevmode
\begin{itemize}
\item \emph{Gradient.} The strong convexity of $r$ gives the strong convexity of $H$ on $K$ with a constant modulus. Therefore, the maximizer is unique and, as in the proof of Lemma~\ref{lem:conj}, $\nabla\psi(y)=z(y)$ and the gradient is continuous. This is statement~(e).
\item \emph{Optimality conditions.} Set $\theta_0=0$ and $\theta_i=y_i+\sigma_i$ for $i\ge1$. The maximization problem in~\eqref{eq:eq-chain} is equivalent to maximizing
\[
 \sum_{i=0}^n\theta_ip_i-\sum_{i=0}^nr(p_i)
\]
over the simplex. Moving weight between two positive coordinates and from a positive coordinate to a zero one gives a number $\nu$ such that for $0\le i\le n$
\begin{equation}\label{eq:eq-kkt}
 \theta_i\ \begin{cases}=\nu+r'(p_i), & p_i>0,\\ \le\nu, & p_i=0.\end{cases}
\end{equation}
These conditions are also sufficient by the convexity of $r$. Since $0\le r'(p_i)\le b$, we always have $\nu\ge\theta_i-b$.
\item \emph{Statement~(a).} The point $z=0$ is feasible, so $\psi\ge0$. Equality holds exactly when the vector $p=e_0$ is optimal. By~\eqref{eq:eq-kkt} this is equivalent to $\nu=-b$ and $\theta_i\le-b$ for $i\ge1$. The formulas for the set of minimizers and the minimizer closest to zero follow from~\eqref{eq:eq-chain-params}.
\item \emph{Statement~(b).} For $y=0$, conditions~\eqref{eq:eq-kkt} hold for $p=e_1$ and $\nu=\sigma_1-b=(2n-2)b$: all other $\theta_i$, including $\theta_0=0$, are at most $\nu$. Therefore, $\psi(0)=G\sigma_1=(4-2/n)\varepsilon$.
\item \emph{Statement~(c).} Substituting $z=Ge_n$, we obtain $\psi(y)\ge G\sigma_n=2\varepsilon$.
\item \emph{Statement~(d).} If $y_i=0$ for $i\ge j$, then for $l>j$
\[
 \theta_l=\sigma_l\le\sigma_j-b=\theta_j-b\le\nu .
\]
A positive $p_l$ would contradict~\eqref{eq:eq-kkt}, since $r'(p_l)>0$. Hence, the $l$-th coordinate of the gradient $Gp_l$ is zero.\qedhere
\end{itemize}
\end{proof}

To verify the class condition, a one-dimensional curvature bound for the function $r$ suffices.

\begin{lemma}[Class condition for equal slopes]\label{lem:eq-conj-class}
If $r$ is twice continuously differentiable on $(0,1)$ and for $0<p<1$
\begin{equation}\label{eq:eq-dual-curvature}
 r''(p)\ge\frac{G}{L_0+L_1Gp},
\end{equation}
then $\psi\in\F_{L_0,L_1}^n$.
\end{lemma}

\begin{proof}
Let $z,z'\in K$ be distinct and $u=z'-z$. On the interval $0<t<1$ the coordinates $p_i$ of the point $z+tu$ are either positive or identically zero. Therefore,
\begin{align*}
 \frac{d^2}{dt^2}H(z+tu)
 &=\frac1G\sum_{i=1}^nr''(p_i)u_i^2+\frac1Gr''(p_0)U_n^2\\
 &\by{\ge}{\eqref{eq:eq-dual-curvature}}\frac{\norm u^2}{L_0+L_1\norm{z+tu}}
 \ge\frac{\norm u^2}{L_0+L_1\max\{\norm z,\norm{z'}\}} ,
\end{align*}
where $U_n=u_1+\dots+u_n$. The terms with identically zero coordinates are omitted. If $p_0\equiv0$, then $U_n=0$ and the last term is absent as well. The second step uses the inequality $Gp_i=(z+tu)_i\le\norm{z+tu}$, and the third one the convexity of the norm. Integration and continuity at the endpoints give the assumption of Lemma~\ref{lem:conj}, which completes the proof.
\end{proof}

\paragraph{Quadratic function $r$.} Take $r(p)=bp^2/2$. Then $r''=b$. Let
\begin{equation}\label{eq:eq-quad-cond}
 4n\le\sqrt{\frac{L_0R^2}\varepsilon} .
\end{equation}
Then
\begin{align*}
 \frac Gb
 &\by{=}{\eqref{eq:eq-chain-params}}\frac{2\varepsilon}{nR^2}\sum_{k=n+1}^{2n}k^2
 \le\frac{8\varepsilon n^2}{R^2}\\
 &\by{\le}{\eqref{eq:eq-quad-cond}}\frac{L_0}2 .
\end{align*}
Consequently, condition~\eqref{eq:eq-dual-curvature} holds, and the gradient of $\psi$ is even Lipschitz with constant at most $L_0/2$.

\paragraph{Logarithmic function $r$.} Set
\begin{equation}\label{eq:eq-log-reg}
\begin{aligned}
 r(p)&=\frac1{L_1}\Bigl[(p+\chi)\log\Bigl(1+\frac p\chi\Bigr)-p\Bigr],\\
 \chi&=\frac1{e^{L_1b}-1} .
\end{aligned}
\end{equation}
Then
\[
\begin{aligned}
 &r'(p)=\frac1{L_1}\log\Bigl(1+\frac p\chi\Bigr),\\
 &r''(p)=\frac1{L_1(p+\chi)},
\end{aligned}
\]
and $r'(1)=b$. Condition~\eqref{eq:eq-dual-curvature} holds if $GL_1\chi\le L_0$, since then $GL_1(p+\chi)\le L_0+L_1Gp$. We verify this under the conditions
\begin{equation}\label{eq:eq-log-cond}
\begin{aligned}
 &L_1^2\varepsilon\ge L_0,\\
 &n\le\Bigl(\frac{L_1R}{8\log(2+L_1^2\varepsilon/L_0)}\Bigr)^{2/3} .
\end{aligned}
\end{equation}
From~\eqref{eq:eq-chain-params} and~\eqref{eq:eq-log-cond}
\begin{align*}
 L_1b
 &\by{\ge}{\eqref{eq:eq-chain-params}}\frac{L_1R}{2n^{3/2}}
 \by{\ge}{\eqref{eq:eq-log-cond}}4\log\Bigl(2+\frac{L_1^2\varepsilon}{L_0}\Bigr)\\
 &\bypad{\ge}{\eqref{eq:eq-log-cond}}{\eqref{eq:eq-chain-params}}4\log3>4 .
\end{align*}
The last step uses the first condition of~\eqref{eq:eq-log-cond}.
Therefore, $e^{L_1b}\ge2$ and
\[
 \chi\le2e^{-L_1b}\le2\Bigl(2+\frac{L_1^2\varepsilon}{L_0}\Bigr)^{-2},
\]
and
\begin{align*}
 \frac{GL_1\chi}{L_0}
 &=\frac{2L_1^2\varepsilon}{nL_1b\,L_0}\,\chi
 \le\frac4{nL_1b}\cdot\frac{L_1^2\varepsilon/L_0}{\bigl(2+L_1^2\varepsilon/L_0\bigr)^2}\\
 &\by{\le}{(\mathrm a)}\frac1{2nL_1b}<\frac12 .
\end{align*}
Here step (a) uses the inequality $t/(2+t)^2\le1/8$ for $t\ge0$, which is equivalent to $(t-2)^2\ge0$.

\begin{theorem}[Lower bound on the class with bounded initial optimality gap]\label{thm:lower-gap}
Let a positive integer $n$ satisfy condition~\eqref{eq:eq-quad-cond} or condition~\eqref{eq:eq-log-cond}. Then for every deterministic method with $n-1$ oracle calls in $\R^{2n}$ there exists a function $f\in\mathcal P_{4\varepsilon}$ on $\R^{2n}$ such that $f^\ast=0$, the distance from zero to the set of minimizers equals $R$, $f(0)=(4-2/n)\varepsilon$ and $f(\hat x)\ge2\varepsilon$. In particular, $\mathcal N(\mathcal P_{4\varepsilon})\ge n$, and for $L_0R^2\ge\varepsilon$
\begin{equation}\label{eq:lower-gap}
 1+\mathcal N(\mathcal P_{4\varepsilon})\ge\frac1{12}\Bigl[1+\sqrt{\frac{L_0R^2}\varepsilon}+\Bigl(\frac{L_1R}{\log(2+L_1^2\varepsilon/L_0)}\Bigr)^{2/3}\Bigr].
\end{equation}
\end{theorem}

\begin{proof}
\leavevmode
\begin{itemize}
\item \emph{Worst-case functions.} If condition~\eqref{eq:eq-quad-cond} holds, take the quadratic function $r$, otherwise the logarithmic one~\eqref{eq:eq-log-reg}. In both cases $\psi\in\F_{L_0,L_1}^n$ by Lemma~\ref{lem:eq-conj-class}, and Lemma~\ref{lem:eq-chain} gives the remaining assumptions of Proposition~\ref{prop:reduction}(b) with $\mathrm{gap}=2\varepsilon$ and $\Delta=4\varepsilon$. By this proposition, for every method there exists a matrix $V$ with orthonormal columns such that the function $f(x)=\psi(V^\top x)$ lies in $\mathcal P_{4\varepsilon}$ and $f(\hat x)\ge2\varepsilon$. The remaining properties of $f$ are given by statements~(a) and~(b) of Lemma~\ref{lem:eq-chain}.
\item \emph{Case $L_1^2\varepsilon\ge L_0$.} For each of the conditions~\eqref{eq:eq-quad-cond} and~\eqref{eq:eq-log-cond}, take as $n$ the integer part of the largest admissible right-hand side. If it is positive, the first part of the theorem gives $\mathcal N(\mathcal P_{4\varepsilon})\ge n$, and if it is zero, this inequality is obvious. Together with the inequality $1+\lfloor x\rfloor\ge x$ we obtain
\begin{align*}
 1+\mathcal N(\mathcal P_{4\varepsilon})&\ge1,\\
 1+\mathcal N(\mathcal P_{4\varepsilon})&\ge\frac14\sqrt{\frac{L_0R^2}\varepsilon},\\
 1+\mathcal N(\mathcal P_{4\varepsilon})&\ge\Bigl(\frac{L_1R}{8\log(2+L_1^2\varepsilon/L_0)}\Bigr)^{2/3}\\
 &=\frac14\Bigl(\frac{L_1R}{\log(2+L_1^2\varepsilon/L_0)}\Bigr)^{2/3} .
\end{align*}
Adding these inequalities with the factors $1$, $4$ and $4$, we obtain that the expression in brackets in~\eqref{eq:lower-gap} is at most $9\bigl(1+\mathcal N(\mathcal P_{4\varepsilon})\bigr)$.
\item \emph{Case $L_1^2\varepsilon<L_0$.} Then $(L_1R)^2<L_0R^2/\varepsilon$, and for $L_0R^2\ge\varepsilon$
\begin{align*}
 \Bigl(\frac{L_1R}{\log(2+L_1^2\varepsilon/L_0)}\Bigr)^{2/3}
 &\by{\le}{(\mathrm a)}(\log2)^{-2/3}(L_1R)^{2/3}
 \by{<}{(\mathrm b)}\frac32\Bigl(\frac{L_0R^2}\varepsilon\Bigr)^{1/3}\\
 &\by{\le}{(\mathrm c)}\frac32\sqrt{\frac{L_0R^2}\varepsilon} ,
\end{align*}
where step (a) uses the inequality $\log(2+t)\ge\log2$, step (b) uses the inequalities $(\log2)^{-2/3}<3/2$ and $(L_1R)^2<L_0R^2/\varepsilon$, and step (c) uses the condition $L_0R^2\ge\varepsilon$. Hence, the expression in brackets in~\eqref{eq:lower-gap} is at most
\begin{align*}
 1+\frac52\sqrt{\frac{L_0R^2}\varepsilon}
 &\le\bigl(1+\mathcal N(\mathcal P_{4\varepsilon})\bigr)+10\bigl(1+\mathcal N(\mathcal P_{4\varepsilon})\bigr)\\
 &=11\bigl(1+\mathcal N(\mathcal P_{4\varepsilon})\bigr).
\end{align*}
Here we used the first two inequalities from the previous case: they were obtained without any condition on the sign of $L_1^2\varepsilon-L_0$.
\end{itemize}
In both cases the constant $1/12$ is valid.
\end{proof}

The exponent $2/3$ arises as follows. The shifts $\sigma_i=(2n-i)b$ form an arithmetic progression with common difference $-b$, so $R^2$ is of order $b^2n^3$. The class condition for the logarithmic function $r$ requires $L_1b$ to be at least a logarithm, as in condition~\eqref{eq:eq-log-cond}. These two bounds give a chain length of order $(L_1R)^{2/3}$ up to a logarithm. The initial optimality gap is less than $4\varepsilon$, and Theorem~\ref{thm:upper-gap}(a) below shows that the exponent $2/3$ is tight up to logarithmic factors. In contrast, the initial optimality gap of the functions of Theorem~\ref{thm:ordered} is not uniformly bounded by $4\varepsilon$ (see the remark after it).

\subsection{Variant of the method without the \texorpdfstring{$\Phi$}{Phi} term}\label{sec:gap-method}

In this subsection and below, $f\in\F_{L_0,L_1}^d$ and $x^\ast$ is a minimizer of $f$. We use a variant of Algorithm~\ref{alg:am} without the additional term $\Phi$ and with an initial weight $a_0\ge0$ (Algorithm~\ref{alg:am-gap}). Set $A_k=a_0+a_1+\dots+a_k$. The coefficient $a_k>0$ is defined by the equation
\begin{equation}\label{eq:coef-gap-rule}
 32\bigl(L_0+L_1\norm{\nabla f(y^k)}\bigr)a_k^2=A_{k-1}+a_k ,
\end{equation}
and the number of line search steps by condition~\eqref{eq:search-m}.

\begin{algorithm}[ht]
\caption{Variant of SHERPA without the $\Phi$ term}\label{alg:am-gap}
\begin{algorithmic}[1]
\Require starting point $x^0$, initial weight $a_0\ge0$, numbers $r,\delta>0$, a positive integer $m$ satisfying~\eqref{eq:search-m}, number of iterations $K$
\State $v^0=x^0$, $A_0=a_0$
\For{$k=1,\dots,K$}
 \State $y^k$ is the output of Algorithm~\ref{alg:search} on the segment $[v^{k-1},x^{k-1}]$ with the numbers $r$, $\delta$, $m$
 \If{the stopping criterion was triggered}
  \State \Return $y^k$
 \EndIf
 \State $a_k$ is the positive root of equation~\eqref{eq:coef-gap-rule}, $A_k=A_{k-1}+a_k$
 \State $x^k=y^k-\dfrac{\nabla f(y^k)}{4(L_0+L_1\norm{\nabla f(y^k)})}$
 \State $v^k=v^{k-1}-a_k\nabla f(y^k)$
\EndFor
\State \Return $x^K$
\end{algorithmic}
\end{algorithm}

\begin{lemma}\label{lem:energy-gap}
Let $\norm{x^0-x^\ast}\le R_0$, let the initial optimality gap $f(x^0)-f^\ast$ be at most $\Delta_0$, and let
\[
 r^2\ge R_0^2+2a_0\Delta_0 .
\]
If Algorithm~\ref{alg:am-gap} completes its first $k$ iterations without meeting the stopping criterion, then
\begin{equation}\label{eq:energy-gap}
 A_k\bigl(f(x^k)-f^\ast\bigr)+\frac12\sum_{i=1}^ka_i^2\norm{\nabla f(y^i)}^2+\frac12\norm{v^k-x^\ast}^2\le\frac12R_0^2+a_0\Delta_0 ,
\end{equation}
and all points of the method and of the line search lie in the ball $B(x^\ast,r)$. Moreover, for all $i\le k$
\begin{equation}\label{eq:gap-monotone}
\begin{aligned}
 f(x^{i-1})&\ge f(y^i),\\
 f(y^i)&\ge f(x^i)+\frac{\norm{\nabla f(y^i)}^2}{8\bigl(L_0+L_1\norm{\nabla f(y^i)}\bigr)} .
\end{aligned}
\end{equation}
\end{lemma}

\begin{proof}
The proof repeats the proof of Lemma~\ref{lem:energy}. For $k=0$, inequality~\eqref{eq:energy-gap} follows from the conditions on $x^0$. The induction step differs in three places.
\begin{itemize}
\item \emph{Replacements.} Instead of Lemma~\ref{lem:anchored}, the convexity inequality $\ip{\nabla f(y^k)}{y^k-x^\ast}\ge f(y^k)-f^\ast$ is used, and Lemma~\ref{lem:search} is used in the same way. Inequalities~\eqref{eq:gap-monotone} follow from it and Lemma~\ref{lem:safe-step}.
\item \emph{Increment.} With these replacements, the bounds on the distance and on the optimality gap from the proof of Lemma~\ref{lem:energy} give
\begin{align*}
 \MoveEqLeft\Bigl[A_k\bigl(f(x^k)-f^\ast\bigr)+\frac12\norm{v^k-x^\ast}^2\Bigr]-\Bigl[A_{k-1}\bigl(f(x^{k-1})-f^\ast\bigr)+\frac12\norm{v^{k-1}-x^\ast}^2\Bigr]\\
 &\bypad{\le}{}{\eqref{eq:coef-gap-rule}}\Bigl(\frac{a_k}{16}-\frac{5A_k}{24}\Bigr)\frac{\norm{\nabla f(y^k)}^2}{L_0+L_1\norm{\nabla f(y^k)}}+\frac12a_k^2\norm{\nabla f(y^k)}^2\\
 &\by{=}{\eqref{eq:coef-gap-rule}}\Bigl(\frac{a_k}{16}-\frac{5A_k}{24}+\frac{A_k}{32}\Bigr)\frac{\norm{\nabla f(y^k)}^2}{L_0+L_1\norm{\nabla f(y^k)}}-\frac12a_k^2\norm{\nabla f(y^k)}^2\\
 &\bypad{\le}{(\mathrm a)}{\eqref{eq:coef-gap-rule}}-\frac12a_k^2\norm{\nabla f(y^k)}^2 ,
\end{align*}
where step (a) uses the inequalities $a_k\le A_k$ and $1/16-5/24+1/32<0$.
\item \emph{Ball of radius $r$.} By inequality~\eqref{eq:energy-gap} and the assumption of the lemma
\[
 \norm{v^k-x^\ast}^2\le R_0^2+2a_0\Delta_0\le r^2 ,
\]
and the point $x^k$ lies in the ball by Lemma~\ref{lem:safe-step}.\qedhere
\end{itemize}
\end{proof}

Equation~\eqref{eq:coef-gap-rule} and the equality $A_k=A_{k-1}+a_k$ imply, for every $i\ge1$, the relations
\begin{equation}\label{eq:coef-gap}
\begin{aligned}
 a_i^2&=\frac{A_i}{32\bigl(L_0+L_1\norm{\nabla f(y^i)}\bigr)},\\
 \sqrt{A_i}-\sqrt{A_{i-1}}&=\frac{a_i}{\sqrt{A_i}+\sqrt{A_{i-1}}}\by{\ge}{A_{i-1}\le A_i}\frac{a_i}{2\sqrt{A_i}}\\
 &=\frac1{2\sqrt{32\bigl(L_0+L_1\norm{\nabla f(y^i)}\bigr)}} .
\end{aligned}
\end{equation}
Moreover, for $A_{i-1}>0$
\begin{equation}\label{eq:coef-gap-log}
 \frac{a_i}{A_i}\le-\log\Bigl(1-\frac{a_i}{A_i}\Bigr)=\log\frac{A_i}{A_{i-1}} .
\end{equation}

\subsection{Optimality gap at most \texorpdfstring{$4\varepsilon$}{4ε}: block bound}\label{sec:count1}

The first counting method bounds the optimality gap after a block of $N$ consecutive iterations in terms of the gap at the beginning of the block.

\begin{lemma}[Block bound]\label{lem:block}
Let $\norm{x^0-x^\ast}\le R$, let Algorithm~\ref{alg:am-gap} be run with $a_0=0$ and $r=R$, and suppose that it completes iterations $k+1,\dots,k+N$, $N\ge1$, without meeting the stopping criterion. Then
\begin{equation}\label{eq:block}
 f(x^{k+N})-f^\ast\le\frac{96L_0R^2}{N^2}+\frac{3072\,L_1^2R^2}{N^3}\bigl(f(x^k)-f^\ast\bigr).
\end{equation}
\end{lemma}

\begin{proof}
All sums in the proof are over $i=k+1,\dots,k+N$.
\begin{itemize}
\item \emph{Optimality gap through the weights.} We have
\begin{align*}
 f(x^{k+N})-f^\ast
 &\by{\le}{\eqref{eq:energy-gap}}\frac{R^2}{2A_{k+N}}
 \le\frac{R^2}{2\bigl(\sqrt{A_{k+N}}-\sqrt{A_k}\bigr)^2}\\
 &\by{\le}{\eqref{eq:coef-gap}}\frac{64R^2}{\Bigl(\sum_i\bigl(L_0+L_1\norm{\nabla f(y^i)}\bigr)^{-1/2}\Bigr)^2} .
\end{align*}
\item \emph{Hölder's inequality.} With exponents $3/2$ and $3$
\begin{align*}
 N&=\sum_i\bigl(L_0+L_1\norm{\nabla f(y^i)}\bigr)^{-1/3}\bigl(L_0+L_1\norm{\nabla f(y^i)}\bigr)^{1/3}\\
 &\le\Bigl(\sum_i\bigl(L_0+L_1\norm{\nabla f(y^i)}\bigr)^{-1/2}\Bigr)^{2/3}\Bigl(\sum_i\bigl(L_0+L_1\norm{\nabla f(y^i)}\bigr)\Bigr)^{1/3} .
\end{align*}
Substituting into the previous bound, we obtain
\begin{equation}\label{eq:block-g}
\begin{aligned}
 f(x^{k+N})-f^\ast
 &\le\frac{64R^2}{N^3}\sum_i\bigl(L_0+L_1\norm{\nabla f(y^i)}\bigr)\\
 &=\frac{64L_0R^2}{N^2}+\frac{64L_1R^2}{N^3}\sum_i\norm{\nabla f(y^i)} .
\end{aligned}
\end{equation}
\item \emph{The sum of gradient norms is small when the optimality gap is small.} Summing inequalities~\eqref{eq:gap-monotone} over the block, we obtain
\[
 \sum_i\frac{\norm{\nabla f(y^i)}^2}{L_0+L_1\norm{\nabla f(y^i)}}\le8\bigl(f(x^k)-f^\ast\bigr).
\]
By the Cauchy--Schwarz inequality
\begin{align*}
 \Bigl(\sum_i\norm{\nabla f(y^i)}\Bigr)^2
 &\le\sum_i\frac{\norm{\nabla f(y^i)}^2}{L_0+L_1\norm{\nabla f(y^i)}}\cdot\sum_i\bigl(L_0+L_1\norm{\nabla f(y^i)}\bigr)\\
 &\le8\bigl(f(x^k)-f^\ast\bigr)\Bigl(NL_0+L_1\sum_i\norm{\nabla f(y^i)}\Bigr).
\end{align*}
If the sum $\sum_i\norm{\nabla f(y^i)}$ exceeded $\sqrt{8NL_0(f(x^k)-f^\ast)}+8L_1(f(x^k)-f^\ast)$, then its square would exceed this quantity times the sum itself, which is at least the right-hand side of the last inequality. Therefore,
\[
 \sum_i\norm{\nabla f(y^i)}\le\sqrt{8NL_0\bigl(f(x^k)-f^\ast\bigr)}+8L_1\bigl(f(x^k)-f^\ast\bigr).
\]
\item \emph{Substitution.} Substitute this bound into~\eqref{eq:block-g}. The second term in~\eqref{eq:block-g} becomes at most
\[
 \frac{64L_1R^2}{N^3}\sqrt{8NL_0\bigl(f(x^k)-f^\ast\bigr)}+\frac{512\,L_1^2R^2\bigl(f(x^k)-f^\ast\bigr)}{N^3} .
\]
By the inequality between the arithmetic and geometric means, the first of these terms is at most
\[
 \frac{32L_0R^2}{N^2}+\frac{256\,L_1^2R^2\bigl(f(x^k)-f^\ast\bigr)}{N^3} .
\]
Therefore,
\[
 f(x^{k+N})-f^\ast\le\frac{96L_0R^2}{N^2}+\frac{768\,L_1^2R^2\bigl(f(x^k)-f^\ast\bigr)}{N^3} .
\]
Since $768\le3072$, we obtain~\eqref{eq:block}.\qedhere
\end{itemize}
\end{proof}

Bound~\eqref{eq:block} explains where the two terms come from: the first one is the usual rate of the accelerated method, and the second one is the price of the nonlinearity of the curvature, proportional to the optimality gap at the beginning of the block. If the initial optimality gap is unknown and may be exponentially large, the second term is useless. This obstacle is circumvented in two ways: by the additional term $\Phi$ if only the radius is known (Section~\ref{sec:upper}), and by a positive initial weight if a bound on the optimality gap is known (next subsection).

\subsection{Known bound on the optimality gap: positive initial weight}\label{sec:count2}

Let $\norm{x^0-x^\ast}\le R$, and let a number $\Delta>\varepsilon$ with $f(x^0)-f^\ast\le\Delta$ be known. Run Algorithm~\ref{alg:am-gap} with the parameters
\begin{equation}\label{eq:weight-params}
\begin{aligned}
 &a_0=\frac{R^2}{2\Delta},\\
 &r=\sqrt2R,\\
 &\delta=\varepsilon .
\end{aligned}
\end{equation}
The assumption of Lemma~\ref{lem:energy-gap} with $R_0=R$ and $\Delta_0=\Delta$ holds, since $R^2+2a_0\Delta=2R^2=r^2$.

\begin{lemma}[Counting with a positive initial weight]\label{lem:count-weight}
If Algorithm~\ref{alg:am-gap} with parameters~\eqref{eq:weight-params} completes $K$ iterations without meeting the stopping criterion and $f(x^K)-f^\ast>\varepsilon$, then
\[
 K<16\sqrt{\frac{L_0R^2}\varepsilon}+21\,(L_1R)^{2/3}\log^{2/3}\frac{2\Delta}\varepsilon .
\]
\end{lemma}

\begin{proof}
The right-hand side of~\eqref{eq:energy-gap} is at most $R^2/2+R^2/2=R^2$. Since $f(x^K)-f^\ast>\varepsilon$, this gives
\begin{equation}\label{eq:count-weight}
\begin{aligned}
 A_K&<\frac{R^2}\varepsilon,\\
 \sum_{i=1}^Ka_i^2\norm{\nabla f(y^i)}^2&\le2R^2 .
\end{aligned}
\end{equation}
We split the iterations into two groups: with a small gradient, $L_1\norm{\nabla f(y^i)}\le L_0$, and with a large one, $L_1\norm{\nabla f(y^i)}>L_0$.
\begin{itemize}
\item \emph{Iterations with a small gradient.} For them $L_0+L_1\norm{\nabla f(y^i)}\le2L_0$, and by~\eqref{eq:coef-gap} each such iteration increases $\sqrt{L_0A_i}$ by at least $1/16$. Hence, their number is at most
\[
 16\sqrt{L_0A_K}\by{<}{\eqref{eq:count-weight}}16\sqrt{\frac{L_0R^2}\varepsilon} .
\]
\item \emph{Iterations with a large gradient.} For them $L_0+L_1\norm{\nabla f(y^i)}<2L_1\norm{\nabla f(y^i)}$, therefore
\begin{align*}
 a_i\norm{\nabla f(y^i)}&\by{=}{\eqref{eq:coef-gap}}\frac{A_i\norm{\nabla f(y^i)}}{32\bigl(L_0+L_1\norm{\nabla f(y^i)}\bigr)a_i}\\
 &\bypad{>}{}{\eqref{eq:coef-gap}}\frac{A_i}{64L_1a_i} .
\end{align*}
This and~\eqref{eq:count-weight} give
\[
 \sum_{L_1\norm{\nabla f(y^i)}>L_0}\Bigl(\frac{A_i}{a_i}\Bigr)^2<(64L_1)^2\cdot2R^2=8192\,L_1^2R^2 .
\]
Since $A_{i-1}\ge a_0>0$, by~\eqref{eq:coef-gap-log}
\begin{align*}
 \sum_{L_1\norm{\nabla f(y^i)}>L_0}\frac{a_i}{A_i}
 &\le\log\frac{A_K}{a_0}
 \by{<}{\eqref{eq:count-weight},\,\eqref{eq:weight-params}}\log\frac{2\Delta}\varepsilon .
\end{align*}
By Hölder's inequality with exponents $3$ and $3/2$, the number of iterations with a large gradient equals
\begin{align*}
 \sum_{L_1\norm{\nabla f(y^i)}>L_0}\Bigl(\frac{A_i}{a_i}\Bigr)^{2/3}\Bigl(\frac{a_i}{A_i}\Bigr)^{2/3}
 &\le\Bigl(\sum\Bigl(\frac{A_i}{a_i}\Bigr)^2\Bigr)^{1/3}\Bigl(\sum\frac{a_i}{A_i}\Bigr)^{2/3}\\
 &<(8192\,L_1^2R^2)^{1/3}\log^{2/3}\frac{2\Delta}\varepsilon\\
 &<21\,(L_1R)^{2/3}\log^{2/3}\frac{2\Delta}\varepsilon .\qedhere
\end{align*}
\end{itemize}
\end{proof}

\subsection{Final upper bounds}\label{sec:gap-upper}

\begin{theorem}[Upper bounds under a bounded initial optimality gap]\label{thm:upper-gap}
Let $f\in\F_{L_0,L_1}^d$, $\norm{x^0-x^\ast}\le R$ and $\varepsilon>0$.
\begin{enumerate}
\item[(a)] If $f(x^0)-f^\ast\le4\varepsilon$, then Algorithm~\ref{alg:am-gap} with $a_0=0$, $r=R$, $\delta=\varepsilon$,
\[
 m=\Bigl\lceil\log_2\Bigl(128\Bigl(1+\frac{L_0R^2}\varepsilon+L_1R\Bigr)^2\Bigr)\Bigr\rceil
\]
and the number of iterations
\begin{align*}
 K&=\max\Bigl\{1,\Bigl\lceil\sqrt{\frac{192L_0R^2}\varepsilon}\,\Bigr\rceil,\bigl\lceil(24576\,L_1^2R^2)^{1/3}\bigr\rceil\Bigr\}\\
 &\le3+14\sqrt{\frac{L_0R^2}\varepsilon}+30(L_1R)^{2/3}
\end{align*}
outputs a point with optimality gap at most $\varepsilon$ after at most $K(m+2)$ oracle calls.
\item[(b)] If a number $\Delta>\varepsilon$ with $f(x^0)-f^\ast\le\Delta$ is known, then Algorithm~\ref{alg:am-gap} with parameters~\eqref{eq:weight-params},
\[
 m=\Bigl\lceil\log_2\Bigl(128\Bigl(1+\frac{2L_0R^2}\varepsilon+\sqrt2L_1R\Bigr)^2\Bigr)\Bigr\rceil
\]
and the number of iterations
\[
 K=\Bigl\lceil16\sqrt{\frac{L_0R^2}\varepsilon}+21\,(L_1R)^{2/3}\log^{2/3}\frac{2\Delta}\varepsilon\Bigr\rceil
\]
outputs such a point after at most $K(m+2)$ oracle calls. For $\Delta\le\varepsilon$ it suffices to return $x^0$.
\end{enumerate}
In both cases the method uses only gradients and does not use the value $f^\ast$. In case~(a) it does not use a bound on the initial optimality gap either.
\end{theorem}

\begin{proof}
\leavevmode
\begin{itemize}
\item \emph{Part~(a).} If the stopping criterion was triggered, the output is correct by Lemma~\ref{lem:search}. Otherwise bound~\eqref{eq:block} with $k=0$ and $f(x^0)-f^\ast\le4\varepsilon$ gives
\[
 f(x^K)-f^\ast\le\frac{96L_0R^2}{K^2}+\frac{3072\,L_1^2R^2\cdot4\varepsilon}{K^3}\by{\le}{(\mathrm a)}\frac\varepsilon2+\frac\varepsilon2 ,
\]
where step (a) uses the choice of $K$. The bound on $K$ follows from the inequalities $\sqrt{192}<14$ and $24576^{1/3}<30$.
\item \emph{Part~(b).} For $\Delta\le\varepsilon$ we return $x^0$. Otherwise, by Lemma~\ref{lem:count-weight}, after $K$ iterations either the stopping criterion was triggered and the output is correct by Lemma~\ref{lem:search}, or $f(x^K)-f^\ast\le\varepsilon$.
\end{itemize}
In both parts the number of calls is at most $K(m+2)$.
\end{proof}

The number of line search steps $m+2$ grows only logarithmically in $L_0R^2/\varepsilon$ and $L_1R$, so the theorem gives, for the class with initial optimality gap at most $4\varepsilon$ and for the class with a known bound $\Delta$ on the initial optimality gap, respectively, the bounds
\begin{align*}
 &O\Bigl(\Bigl[1+\sqrt{\frac{L_0R^2}\varepsilon}+(L_1R)^{2/3}\Bigr]\log\Bigl(2+\frac{L_0R^2}\varepsilon+L_1R\Bigr)\Bigr),\\
 &O\Bigl(\Bigl[1+\sqrt{\frac{L_0R^2}\varepsilon}+(L_1R)^{2/3}\log^{2/3}\frac{2\Delta}\varepsilon\Bigr]\log\Bigl(2+\frac{L_0R^2}\varepsilon+L_1R\Bigr)\Bigr).
\end{align*}
In case~(b), before the run, $\Delta$ may be replaced by the minimum of $\Delta$ and the right-hand side of~\eqref{eq:radial-gap}. Together with Theorem~\ref{thm:lower-gap}, the first bound gives for $L_0R^2\ge\varepsilon$
\[
 1+\mathcal N(\mathcal P_{4\varepsilon})=\widetilde\Theta\Bigl(1+\sqrt{\frac{L_0R^2}\varepsilon}+(L_1R)^{2/3}\Bigr).
\]
A lower bound matching~(b) for $\Delta$ between $4\varepsilon$ and the right-hand side of~\eqref{eq:radial-gap} is unknown.

\section{Fixed dimension}\label{app:fixed-dim}

The lower bounds of Section~\ref{sec:lower} require a dimension that grows with the number of oracle calls. We show that this is essential: in fixed dimension, the center-of-gravity method solves the problem with fixed normalized accuracy $L_1^2\varepsilon/L_0$ in a number of calls logarithmic in $L_1R$.

The function $t\mapsto e^t-1-t$ strictly increases on $[0,\infty)$ from $0$ to $+\infty$. Therefore, the number $\rho_\varepsilon>0$ is uniquely defined as the radius of the ball around the minimizer in which the growth bound of Lemma~\ref{lem:radial} does not exceed $\varepsilon$:
\begin{equation}\label{eq:fixed-radius}
 \frac{L_0}{L_1^2}\bigl(e^{L_1\rho_\varepsilon}-1-L_1\rho_\varepsilon\bigr)=\varepsilon .
\end{equation}

\begin{proposition}[Center-of-gravity method]\label{prop:fixed-dim}
Let $L_0,L_1,R,\varepsilon>0$ and let the dimension $d$ be fixed. The class $\mathcal P_{\rm R}$ in dimension $d$ is solved by the center-of-gravity method with
\begin{equation}\label{eq:centroid-count}
 1+\Bigl\lfloor\frac{d\log\bigl(2\sqrt d\,R/\min\{R,\rho_\varepsilon\}\bigr)}{-\log(1-1/e)}\Bigr\rfloor
\end{equation}
oracle calls, that is, with
\[
 O\Bigl(d\log\Bigl(2\sqrt d\max\Bigl\{1,\frac R{\rho_\varepsilon}\Bigr\}\Bigr)\Bigr)
\]
calls. In dimension $d=1$, bisection uses at most
\begin{equation}\label{eq:one-dim-count}
 \Bigl\lceil\max\Bigl\{\log_2\frac R{\rho_\varepsilon},\,0\Bigr\}\Bigr\rceil
\end{equation}
calls. If the output must coincide with a query point, one additional call suffices.
\end{proposition}

For fixed $d$ and $L_1^2\varepsilon/L_0$, the quantity $L_1\rho_\varepsilon$ is constant, and bound~\eqref{eq:centroid-count} is of order $\log(2+L_1R)$. The multidimensional center-of-gravity method selects the output by function values. Bisection uses only the signs of the derivative.

\begin{proof}
\leavevmode
\begin{itemize}
\item \emph{Ball of points with small optimality gap.} By Lemma~\ref{lem:radial}, for every $x$
\begin{equation}\label{eq:fixed-envelope}
 f(x)-f^\ast\le\frac{L_0}{L_1^2}\bigl(e^{L_1\norm{x-x^\ast}}-1-L_1\norm{x-x^\ast}\bigr).
\end{equation}
Consequently, in the ball $B(x^\ast,\min\{R,\rho_\varepsilon\})$ the optimality gap is at most $\varepsilon$.
\item \emph{Center-of-gravity method.} Without loss of generality $x^0=0$. Start with the cube $K_0=[-2R,2R]^d$. It contains the ball $B(x^\ast,\min\{R,\rho_\varepsilon\})$, since $\norm{x^\ast}\le R$. At step $k$, compute the center of gravity $c_k$ of the polytope $K_{k-1}$ and call the oracle at this point. If $\nabla f(c_k)=0$, then $c_k$ is a minimizer. Otherwise set
\[
 K_k=K_{k-1}\cap\bigl\{z:\ \ip{\nabla f(c_k)}{z-c_k}\le0\bigr\}.
\]
After the prescribed number of steps, return the query point with the smallest function value. Computing the centers of gravity between oracle calls is not counted.
\item \emph{The ball is preserved.} Suppose that at all of the first $N$ query points the optimality gap exceeds $\varepsilon$. For every point $z$ of the ball $B(x^\ast,\min\{R,\rho_\varepsilon\})$, convexity and~\eqref{eq:fixed-envelope} give
\[
 \ip{\nabla f(c_k)}{z-c_k}\le f(z)-f(c_k)<0 .
\]
Hence, the whole ball is preserved after each cut.
\item \emph{Comparison of volumes.} By Grünbaum's theorem \citep{grunbaum1960partitions}, every hyperplane through the center of gravity of a convex body leaves in each of the two half-spaces at least $1/e$ of its volume. The unit ball contains the cube $[-1/\sqrt d,1/\sqrt d]^d$, so its volume is at least $(2/\sqrt d)^d$. Consequently,
\begin{align*}
 \Bigl(\frac{2\min\{R,\rho_\varepsilon\}}{\sqrt d}\Bigr)^d
 &\by{\le}{(\mathrm a)}\operatorname{vol}(K_N)\\
 &\by{\le}{(\mathrm b)}\Bigl(1-\frac1e\Bigr)^N(4R)^d ,
\end{align*}
where step (a) uses the preservation of the ball and step (b) uses Grünbaum's theorem.
For $N$ equal to~\eqref{eq:centroid-count}, the right-hand side is strictly smaller than the left-hand side. This contradiction shows that at least one of the query points has optimality gap at most $\varepsilon$.
\item \emph{One-dimensional case.} Start with the segment $[-R,R]$. At each step, call the oracle at the midpoint of the current segment. If the derivative there is positive, keep the left half, if it is negative, keep the right half, and if it is zero, an exact minimizer has been found. The monotonicity of the derivative of a convex function guarantees that the minimizer stays in the segment. After $n$ calls, the midpoint $\hat x$ of the remaining segment satisfies
\begin{align*}
 |\hat x-x^\ast|&\bypad{\le}{}{\eqref{eq:fixed-envelope}}R\,2^{-n},\\
 f(\hat x)-f^\ast&\by{\le}{\eqref{eq:fixed-envelope}}\frac{L_0}{L_1^2}\bigl(e^{L_1R2^{-n}}-1-L_1R\,2^{-n}\bigr).
\end{align*}
For $n$ equal to~\eqref{eq:one-dim-count}, we have $R2^{-n}\le\rho_\varepsilon$, and the right-hand side is at most $\varepsilon$. The oracle may not have been called at the output midpoint. If the output must coincide with a query point, one more call is needed.\qedhere
\end{itemize}
\end{proof}

\begin{remark}\label{rem:fixed-dim-growth}
For fixed $d$ and $L_1^2\varepsilon/L_0$, Proposition~\ref{prop:fixed-dim} gives only $O(\log(2+L_1R))$ calls. Therefore, a lower bound that is linear in $L_1R$ up to logarithms cannot be realized by a family of fixed dimension. It requires a mechanism that hides information in a number of directions growing together with the number of calls. This is exactly why Theorems~\ref{thm:lower-gap} and~\ref{thm:lower-radius} use a dimension proportional to the chain length.
\end{remark}

\clearpage
\section{Summary table of known bounds}\label{app:literature}
\relpenalty=500 \binoppenalty=700

Table~\ref{tab:lit-convex} collects known upper bounds for the convex deterministic case. Here $R=\norm{x^0-x^\ast}$, $F_0=f(x^0)-f^\ast$, the notation $O(\cdot)$ hides absolute constants, and $\widetilde O(\cdot)$ also hides logarithmic factors. In the ``Smooth.'' column, sym denotes symmetric $(L_0,L_1)$-smoothness (Definition~\ref{def:class}), asym denotes local asymmetric smoothness, $\ell$ denotes smoothness with respect to a general function $\ell(\norm{\nabla f})$, and $(H_0,H_1)$ denotes the condition $\norm{\nabla^2f}\le H_0+H_1(f-f^\ast)$.

\begin{center}
\footnotesize
\captionof{table}{Convex deterministic case: number of oracle calls until $f(\hat x)-f^\ast\le\varepsilon$.}
\label{tab:lit-convex}
\begin{tabularx}{\linewidth}{@{}P{.2\linewidth}lP{.2\linewidth}Y@{}}
\toprule
Paper & Smooth. & Method & Complexity\\
\midrule
\citet{koloskova2023revisiting} & sym $+$ $L$-smooth & clipped GD & $O\bigl(\frac{L_0R^2}{\varepsilon}+L_1R^2\sqrt{\frac L\varepsilon}\bigr)^{a}$\\
\citet{takezawa2024polyak} & sym $+$ $L$-smooth & GD with Polyak stepsizes & $O\bigl(\frac{L_0R^2}{\varepsilon}+L_1R^2\sqrt{\frac L\varepsilon}\bigr)^{a}$\\
\citet{li2023generalized} & $\ell$ & GD & $O\bigl(\frac{\ell(G)R^2}{\varepsilon}\bigr)$, $G=\norm{\nabla f(x^0)}^{f}$\\
\citet{li2023generalized} & $\ell$ & Nesterov's accelerated gradient method & $O\bigl((1+L_1^2R^2+\frac{L_1^2F_0}{L_0})\sqrt{\frac{F_0+L_0R^2}\varepsilon}\bigr)^{g}$\\
\citet{gorbunov2024methods} & sym & $(L_0,L_1)$-GD, GD with Polyak stepsizes & $O\bigl(\frac{L_0R^2}{\varepsilon}+L_1^2R^2\bigr)$\\
\citet{gorbunov2024methods} & sym & $(L_0,L_1)$-STM & $O\bigl(\sqrt{(1+L_1Re^{L_1R})\frac{L_0R^2}{\varepsilon}}\bigr)$\\
\citet{vankov2025optimizing} & sym & GD (three stepsize rules), clipping & $O\bigl(\frac{L_0R^2}{\varepsilon}+L_1R\log\frac{F_0}{\varepsilon}\bigr)$\\
\citet{vankov2025optimizing} & sym & normalized GD, GD with Polyak stepsizes & $O\bigl(\frac{L_0R^2}{\varepsilon}+L_1^2R^2\bigr)$\\
\citet{vankov2025optimizing} & sym & AGMsDR$^{h}$ & $O\bigl(\sqrt{\frac{L_0R^2}{\varepsilon}}+(L_1R)^{2/3}\log\frac{F_0}{\varepsilon}\bigr)$ iterations\\
\citet{lobanov2025linear} & asym & GD, clipped GD & $O\bigl(\frac{L_0R^2}{\varepsilon}+L_1R\log\frac{F_0}{\varepsilon}\bigr)$\\
\citet{lobanov2025linear} & asym & normalized GD & $O\bigl(\bigl(\frac{L_0R^2}{\varepsilon}+L_1R\bigr)\log\frac{F_0}{\varepsilon}\bigr)$\\
\citet{tyurin2025unified} & $\ell$ & GD with stepsize $\int_0^1\frac{dv}{\ell(\norm{\nabla f(x^k)}(1+v))}$ & $O\bigl(\frac{L_0R^2}{\varepsilon}+\min\bigl\{L_1R\sqrt{\frac{F_0}\varepsilon},L_1^2R^2,\frac{L_1\norm{\nabla f(x^0)}R^2}{\varepsilon}\bigr\}\bigr)$\\
\citet{krivchenko2026computer} & sym & GD, PEP & refined constants in the bound with $L_1R\log\frac{F_0}{\varepsilon}$\\
\citet{oikonomidis2025nonlinearly} & anisotropic & nonlinearly preconditioned GD & $O(\mathcal D_0/\varepsilon)$, $\mathcal D_0$ is a supremum over the sublevel set\\
\citet{kovalev2026normalized} & sym & GD with Polyak stepsizes and inexact normalized gradients & $O\bigl(\frac{ML_0R^2}{\varepsilon}+\frac{ML_1^2R^2}{c}\bigr)^{k}$\\
\citet{tyurin2026accelerated} & $\ell$ & accelerated method & $O\bigl(\sqrt{\frac{L_0R^2}{\varepsilon}}+L_1R\log_+\min\bigl\{\frac{L_1^2F_0}{L_0},\frac{F_0}{\varepsilon}\bigr\}\bigr)^{i}$\\
\citet{borodich2026graal} & sym & accelerated GRAAL & $O\bigl(1+\sqrt{\frac{L_0R^2}{\varepsilon}}+(L_1R)^3\bigr)^{j}$\\
\citet{lobanov2026accelerated} & $(H_0,H_1)$ & accelerated, low-dimensional relaxation, restarts & $\widetilde O\bigl(\sqrt{\frac{H_0\tilde R^2}{\varepsilon}}+\sqrt{H_1\tilde R^2}\log\frac{F_0}{\varepsilon}\bigr)^{l}$\\
\midrule
\textbf{This paper} (Theorem~\ref{thm:upper}) & sym & accelerated method with inexact line search & $5\sqrt{\frac{L_0R^2}{\varepsilon}}+O\bigl((1+L_1R)\log(2+L_1R)\bigr)$\\
\bottomrule
\end{tabularx}
\par\smallskip\raggedright
$a$: derived by us from the theorem ($c=L_0/L_1$, $L$-smoothness is additionally assumed). $f$: under the radius condition from the definition of $(r,\ell)$-smoothness. $g$: as recomputed by \citet{vankov2025optimizing}. $h$: each iteration contains an exact line search. $i$: plus an initialization-dependent term. $j$: with the initial stepsize $\eta_0=e^{-L_1R}/L_0$, the full bound also contains $(1+(L_1R)^2)\log(1/(\eta_0L_0))$. $k$: $M$ bounds the gradient norm along the trajectory, $0<c<\min\{\varepsilon/R,L_0/L_1\}$, and the gradient error is small enough. $l$: $\tilde R$ is the radius of the initial sublevel set around $x^\ast$.
\end{center}

\end{document}